\documentclass[leqno,11pt, a4]{amsart}
\usepackage[utf8]{inputenc}
\usepackage{amsthm}
\usepackage{amsmath}
\usepackage{color}
\usepackage{mathtools}
\usepackage{enumitem}
\usepackage[colorinlistoftodos,size=prescript{}]{todonotes}
\usepackage{amsfonts}
\usepackage{amssymb}
\usepackage{hyperref}
\usepackage{pdfsync}
\usepackage{graphicx}
\usepackage{comment}
\usepackage{thmtools}
\usepackage{thm-restate}
\usepackage{cleveref}

\def\R{\text{$\mathbb{R}$}}

\def\fis#1{\dot{#1}}
\def\lra{\longrightarrow}

\def\fis#1{\dot{#1}}
\def\lra{\longrightarrow}

\def\d1#1#2{\frac{d#1}{d#2}}

\def\p1#1#2{\frac{\partial #1}{\partial #2}}

\def\C{\text{$\mathbb{C}$}}
\def\N{\text{$\mathbb{N}$}}

\newcommand\mf[1]{\mathfrak{#1}}

\newcommand\lie[1]{\mathcal L_{#1}}
\newcommand\contr{\it{i}}
\newcommand\lieg{\mathfrak g}

\newcommand\mrm{\mathrm}
\newcommand{\core}{{}_cM}

\newcommand\dist[2]{\mrm{d}_{M/G}({#1},{#2})}
\def\lra{\longrightarrow}
\title[]{Reduction principles for proper actions}
\author{Leonardo Biliotti, Gustavo May Custodio and Alessandro Minuzzo }
\address{Dipartimento di Scienze Matematiche, Fisiche e Informatiche \\
          Universit\`a di Parma (Italy)}
\email{leonardo.biliotti@unipr.it}
\email{victorgustavo.maycustodio@unipr.it}
\email{alessandro.minuzzo@unipr.it}
\thanks{The authors were partially supported by ``National Group for Algebraic and Geometric Structures, and their Applications'' (GNSAGA - INDAM). The first author was partially supported by PRIN  2017
   ``Real and Complex Manifolds: Topology, Geometry and holomorphic dynamics ''.}
\begin{document}

\newtheorem{red}{Theorem}
\renewcommand*{\thered}{\Alph{red}}

\newtheorem{thm}{Theorem}[section]

\newtheorem*{redg}{Theorem (Reduction Principle for Proper Actions)}
\newtheorem*{redc}{Theorem [Reduction Principle for Coisotropic actions]}
\newtheorem*{redalmost}{Theorem [Reduction Principle for infinitesimally almost homogeneous actions]}
\newtheorem*{eqt}{Theorem [Equivalence Theorem for Hamiltonian actions]}
\newtheorem*{equivalence-nonHamiltonian}{Theorem [Equivalence Theorem]}
\newtheorem*{cauzzi}{Restriction Lemma}
\newtheorem{prop}[thm]{Proposition}
\newtheorem{lemma}[thm]{Lemma}
\newtheorem{cor}[thm]{Corollary}
\theoremstyle{definition}
\newtheorem{defini}[thm]{Definition}
\newtheorem{notation}[thm]{Notation}
\newtheorem{exe}[thm]{Example}
\newtheorem{conj}[thm]{Conjecture}
\newtheorem{prob}[thm]{Problem}
\theoremstyle{remark}
\newtheorem{rem}[thm]{Remark}
\newcommand{\thistheoremname}{}
\theoremstyle{plain}
\newtheorem*{genericthm*}{\thistheoremname}
\newenvironment{namedthm*}[1]
  {\renewcommand{\thistheoremname}{#1}%
   \begin{genericthm*}}
  {\end{genericthm*}}
\newcommand{\restr}[1]{\vert_{#1}}
\newcommand{\act}[2]{#1 \circlearrowright #2}
\newcommand{\actg}{G \circlearrowright X}
\newcommand{\actgl}{G \circlearrowright M}
\newcommand{\orbg}{G(p)}
\newcommand{\euclideo}{\langle \cdot , \cdot \rangle}
\def\dq{\mathrm{d}_{M/G}}
\def\dgr{\mathrm{d}_{M^H/N(H)}}
\newcommand{\princ}{M_{\mathrm{princ}}}
\begin{abstract}
We study the \emph{core} of a proper action by a Lie group $G$ on a smooth manifold $M$, extending the construction for $G$ compact by Skjelbred and Straume in \cite{SS}.
Moreover, we show that many properties of a proper $G$-action on $M$ are determined by the action of a group $G'$ on the corresponding core $\core$. We say that such properties admit a \emph{reduction principle}. 
In particular, we prove that a proper isometric $G$-action on $M$ is polar (resp. hyperpolar) if and only if the $G'$-action on $\core$ is polar (resp. hyperpolar).  
In the case of a proper action by syplectomorphisms on a symplectic manifold, we show that a reduction principle holds for coisotropic and infinitesimally almost homogeneous actions.  
We further study the coisotropic condition for the case of a proper Hamiltonian action and its relation with the symplectic stratification described by Lermann and Bates in \cite{lb}.
In particular, we obtain several characterizations for coisotropic actions, some of which extend known results for the action of a compact group of holomorphic automorphisms on a compact K\"ahler manifold obtained by Huckleberry and Wurzbacher in \cite{HW}.
Finally, we study some applications of the core construction for the action of a compact Lie group on a K\"ahler manifold by holomorphic isometries.

\end{abstract}
\maketitle
\section{Introduction}
Let $G$ be a Lie group acting on a smooth manifold $M$. The manifold $M$, together with a given $G$-action, is also called a \emph{$G$-manifold} which we denote with the pair $(G,M)$. We will always assume that the $G$-action is proper. Then the stratification of $M$ by orbit types induces a Whitney stratification on the quotient $M/G$, see \cite{BM} and \cite{S80}. In the literature  (for instance, in \cite{got,pg,Gorodski,Path1} and \cite{SS}) the term \emph{reduction} of a $G$-manifold $M$ is used to indicate a $G'$-manifold, $(G',M')$, where $G'$ is a Lie group acting on a immersed submanifold $M'\subseteq M$ such that the immersion $M'\hookrightarrow M$ induces an isomorphism $M'/G'\lra M/G$ of stratified spaces. When $G$ acts by isometries on a Riemannian manifold $(M,\mathtt g)$, where $\mathtt g$ is the Riemannian metric, the quotient $M/G$ inherits a metric space structure for which the quotient map $\pi:M\lra M/G$ is a submetry (See Section \ref{sect preliminaries}). In this case, we refer to the reduction $(G',M')$ of $(G,M)$ as a \emph{metric reduction} if $G'$ acts on $(M',\mathtt g\vert_{M'})$ by isometries and the immersion $(M',\mathtt g\vert_{M'})\hookrightarrow (M,\mathtt g)$ induces an isometry $M'/G'\lra M/G$.

A classical example is the reduction of the action of a compact connected Lie group $K$ on itself by conjugation to a maximal torus and its associated Weyl group, see \cite{BM} and \cite{DK}. This is a particular case of a \emph{polar action} (see Section \ref{sect polar}). It is well known (see e.g. \cite{BM} or \cite{Gorodski}) that a polar $G$-manifold $M$ admits a metric reduction $(W,\Sigma)$, where $W$ is a discrete group. In this case, the dimension of $\Sigma$ equals the dimension of $M/G$. 



In the unpublished paper \cite{SS}, Skjelbred and Straume considered the action of a compact Lie group $G$ on a connected smooth manifold $M$ and constructed a submanifold $\core\subseteq M$, called the \emph{core} of $(G,M)$, and a group $L$ for which $(L,\core)$ is a reduction of $(G,M)$. In \cite{GS}, Grove and Searle showed that for any $G$-invariant metric on $M$, $(L,\core)$ is also a metric reduction and used it to study the geometry of $M$. We will call the couple $(L,\core)$ the \emph{reduction to the core} or the \emph{core reduction} of $(G, M)$.  We also remark that Luna and Richardson studied in \cite{LR} the core reduction in the algebraic setting. In this paper, we study the reduction to the core in the more general setting of proper actions.

We will always assume that the quotient $M/G$ is connected. Then, there exists only one principal orbit type; see \cite[Theorem 2.8.5]{DK}. Throughout the paper, we let $H$ be a principal isotropy and define $\core$ to be the closure of $M^H \cap \princ$ in $M$, where $M^H$ is the fixed point set of $H$ and $\princ$ is the union of all the principal orbits. The set $\core$ is a closed, embedded submanifold of $M$ (see Remark \ref{core i union of conn.comps}); following \cite{SS} and \cite{GS} we call it the \emph{core} of $(G,M)$. Let $N$ denote the normalizer $N(H)$ of $H$ in $G$. The $G$-action on $M$ induces an $N$-action on $\core$. Since $H$ acts trivially on $\core$, then $L\coloneqq N/H$ acts effectively on $\core$. In Section \ref{sect reduction to core} we prove the following

\begin{restatable}[Reduction to the Core for Proper Actions]{red}{corereduction}\label{stratifiediso}
    The inclusion $\core\hookrightarrow M$ induces an isomorphism  $\Theta:\core/L \lra M/G$ of stratified spaces. In particular, $(\core)_{\mrm{princ}}=M_H$.
\end{restatable}
This, together with Theorem \ref{reduction-global}, implies that for any $G$-invariant metric on $M$, $(L,\core)$ is a metric reduction of $(G,M)$. For the sake of completeness we also provide the proof of results stated without proof in \cite{SS,Straume2}.

\begin{rem}
    It is well known that the connected component of the identity, $N^0$, of $N$ coincides with the connected component of the identity of $Z(H)H$, where $Z(H)$ is the centralizer of $H$ in $G$. It follows that the connected components of the $L$-orbits in $\core$ are the same as the $Z(H)^0$-orbits in $\core$.
\end{rem}
 Furthermore, let $F$ be a connected component of $\core$ and let $L_F$ be the maximal subgroup of $L$ that preserves $F$. We show in Theorem \ref{reduction2} that also $(F,L_F)$ is a metric reduction of $(G,M)$.
 
 Let $C^\infty(M)^G$ denote the algebra of smooth functions on $M$ that are constant along the $G$-orbits in $M$. In Corollary \ref{ChevCore} we show that the restriction map $C^\infty (M)\lra C^\infty (\core)$ induces an isomorphism $ C^\infty (M)^G\lra C^\infty (\core)^L$ of algebras. Corollary \ref{ChevCore} exemplifies one aspect of the classical process of reduction of variables, for which $G$-invariant objects on $M$ are in bijection with $L$-invariant objects on $\core$. 
 
 One of the most useful consequences of the reduction to the core is the so-called \emph{reduction principle}. This involves general properties of the $G$-action on a manifold $M$, and how these properties descend to the $L$-action on the core. Formally, let $\mathcal P$ be a property for $G$-spaces (e.g. a $G$-manifold being polar, of cohomegeneity $k$, or coisotropic, etc.). We say that $\mathcal P$ admits a reduction principle if the following holds
\begin{quotation}
     \emph{The $G$-action on $M$ satisfies property $\mathcal P$ if and only if the $L$-action on $\core$ satisfies property $\mathcal P$.}
\end{quotation}
In this paper, we exhibit some explicit instances of the reduction principle for proper actions. We start with polar and hyperpolar proper actions.

Recall that a $G$-action on $M$ is called \emph{polar} if there exists a complete connected immersed submanifold $\Sigma$ of $M$ which intersects all the orbits and is orthogonal to the $G$-orbit in all common points. Such submanifold $\Sigma$ is automatically totally geodesic and it is called a \emph{section}. If there is a section which is flat in its induced Riemannian metric, then the action is called
\emph{hyperpolar}. Polar and hyperpolar actions have been studied by Conlon \cite{conlon} and Palais and Terng \cite{Path1} among many others. Natural examples for hyperpolar actions are given by the isotropy actions and isotropy representations of symmetric spaces.
In Section \ref{sect polar} we prove that both the polar and the hyperpolar conditions admit a reduction principle.
\begin{restatable}[Reduction Principle for Polar Actions]{red}{polarred}\label{reduction-polar}
Let $G$ be a Lie group acting properly on $M$.
Then the $G$-action on $M$ is polar (hyperpolar, resp.) if and only if the $L$-action on $\core$ is polar (hyperpolar, resp.).
\end{restatable} 
Theorem \ref{reduction-polar} is well known to experts; it is stated, e.g., in \cite{GZ}; however, to the best of our knowledge, no proof is available in the literature.  Theorem \ref{reduction-polar} and the fact that the slice representation at some $x\in \core$ for $(L,\core)$ is the reduction to the core of the slice representation at $x$ for $(G,M)$ (see Section \ref{sect reduction to core}) imply Corollary \ref{infpolar-reduction} which states that the property of every slice representation being polar admits a reduction principle. 


From Section \ref{sa} on we consider the proper action of a Lie group $G$ by symplectomorphisms on a symplectic manifold $(M,\omega)$, where $\omega$ is the symplectic form. The $G$-action on $M$ is called \textit{coisotropic} if there exists an open dense subset $U$ of $M$ such that $G(y)$ is a coisotropic submanifold of $M$ for all $y\in U$.
Coisotropic actions are studied in \cite{multifree} and \cite{Kac} with the terminology \emph{multiplicity free}, which is justified by certain representation theoretic applications. Another motivation for the study of coisotropic actions comes from K\"ahler Geometry. Assume that  $(M,\omega)$ is a compact connected K\"ahler manifold and $G$ is a compact connected Lie group acting  isometrically on $M$. This action is automatically holomorphic by a theorem of Kostant \cite{KN} and it induces, by compactness of $M$, an action on $M$ of the complexified group $G^\C$, in the sense of \cite{ho}. If the $G$-action on $M$ is Hamiltonian (see Section \ref{sect hamiltonian}), then the $G$-action is coisotropic if and only if every Borel subgroup of the complexified group $G^\C$  has an open orbit in $M$. Such open $G^\C$-orbit $\Omega$ is called a \emph{spherical homogeneous space} and $M$ is called a \emph{spherical embedding} of $\Omega$, see for instance  \cite{ak,brion,brion2} and \cite{HW}. Moreover, by \cite[Equivalence Theorem p.274]{HW}, the K\"ahler manifold is projective algebraic. We show that the coisotropic condition admits a reduction principle.  

\begin{restatable}[Reduction Principle for Coisotropic Actions]{red}{coisotropicreduction} \label{reduction-principle-coiso}

Let $G$ be a Lie group acting properly on $(M,\omega)$ by symplectomorphisms.
Then the $G$-action on $M$ is coisotropic if and only if the $L$-action on $\core$ is coisotropic.
    
\end{restatable}

Since $G$ acts by symplectomorphism, it is easy to prove that the action preserves the Poisson bracket, in the sense that for all $f,h\in C^\infty(M)$ and $g\in G$, $g^*\{f,h\}=\{g^*f,g^*h\}$. Therefore, the set of $G$-invariant smooth functions $C^\infty (M)^G$ is a Poisson subalgebra of $\left ( C^\infty (M),\{\cdot,\cdot \} \right )$. We prove in Proposition \ref{multfree iff coiso} that the $G$-action on $M$ is coisotropic if and only if $C^{\infty} (M)^G$ is Abelian. This result is known  under the assumption that $G$ is compact and the $G$-action on $M$ is Hamiltonian \cite{Ben,DP2,multifree} and \cite{HW}. We use this to obtain the main result of Section \ref{sa}.

\begin{restatable}[Equivalence Theorem for Coisotropic Actions]{red}{equivalencethmcoisotropic}\label{equivalence}
Let $G$ be a Lie group acting properly on $(M,\omega)$ by symplectomorphisms. Let $G^0$ and $L^0$ be the identity components of $G$ and $L$, respectively.
 Then the following conditions are equivalent:
\begin{enumerate}
\item $G$ acts coisotropically on $M$;
\item $G^o$ acts coisotropically on a connected component of $M$;
\item for any $x\in \princ$, the orbit $G(x)$ is coisotropic;
\item $C^\infty (M)^G$ is an Abelian Lie algebra;
\item $L$ acts coisotropically on $\core$;
\item $L^o$ acts coisotropically on a connected component of $\core$.
\end{enumerate}
    
\end{restatable}
As noted above, if $G$ is compact and $M$ is a compact K\"ahler manifold, then $G^\mathbb C$ acts on $M$ holomorphically. If, moreover, the $G$-action on $M$ is coisotropic, then there is an open $G^\mathbb C$-orbit, and the $G^\C$-action on $M$ is called \emph{almost homogeneous}. Since in our setting $M$ is not necessarily compact, in general, the complexified group $G^\C$ does not act on $M$. Nevertheless, one can introduce a weaker notion of almost homogeneity as follows. We say that the $G$-action on the almost K\"ahler manifold $(M,\omega,J)$ is \emph{infinitesimally almost homogeneous} if there exists some $x\in M$ for which $T_xM=T_xG(x) +J\left( T_xG(x)\right)$. We show in Section \ref{sect inf almost hom} that this property admits a reduction principle.

\begin{restatable}[Reduction Principle for Infinitesimally Almost Homogeneous Actions]{red}{infhomreduction}\label{red infhom} Let $G$ be a Lie group acting properly on $(M,\omega)$ by symplectomorphisms.
The $G$-action on $M$ is infinitesimally almost homogeneous if and only if the $L$-action on $\core$ is infinitesimally almost homogeneous.
\end{restatable}

In Section \ref{sect hamiltonian} we study the coisotropy condition under the additional hypotheses that the $G$-action on $(M,\omega,J)$ is Hamiltonian with a $G$-equivariant momentum map $\mu:M\lra\mathfrak g^*$. We use the local description of the momentum map and the stratification of the Hamiltonian reductions developed by Bates and Lerman in \cite{lb} to relate the coisotropic condition in the Hamiltonian setting to the Hamiltonian reduction. In particular, we recall that given $\alpha\in\mathfrak{g}^*$, the so-called \emph{reduced space} $M_\alpha:=\mu^{-1}(G(\alpha))/G$, where $G(\alpha)\subseteq \mathfrak{g}^*$ is the coadjoint orbit through $\alpha$, is a stratified space whose strata are symplectic manifolds. We can now state the main result of Section \ref{sect hamiltonian}.

\begin{restatable}[Equivalence Theorem for Coisotropic Hamiltonian Actions]{red}{hamiltoniancoisoreduction}
    \label{equiham}
    Let $(M,\omega)$ be a symplectic manifold on which a Lie group $G$ acts properly and in a Hamiltonian fashion, with momentum map $\mu:M\to \mathfrak g^*$. Assume also that the coadjoint orbits of $G$ are locally closed and $M/G$ is connected. Then the following conditions are equivalent. 

    \begin{itemize}
        \item[(1)-(6)] of Theorem \ref{equivalence};
        \item[(7)] for every $\alpha\in \mathfrak g^*$, the reduced space $M_\alpha$ is discrete;
        \item[(8)] for any $x\in \princ$, we have that 
        \begin{align*}
            \dim (G_{\mu(x)})- \dim (G_x)=\text{cohom}(G,M);
        \end{align*}
        \item[(9)] let $\mathfrak l$ be the Lie algebra of $L$. For every $\alpha\in \mathfrak l^*$, the reduced space $(\core)_\alpha$ is discrete.
    \end{itemize}
\end{restatable}

Finally, in Section \ref{sect kahler} we prove a reduction principle for Hamiltonian actions on a connected K\"ahler manifold. 

Let $(M,\omega, J,\mathtt h)$ be a K\"ahler manifold and let $G$ be a compact connected Lie group. We assume that the $G$-action on $M$ is Hamiltonian with a $G$-equivariant momentum map $\mu:M\lra \lieg^*$. We also assume that the $G$-action is by holomorphich isometries and the $G$-action extends holomorphically to an action of the complexified group $G^\C$. 
We fix a $G$-invariant scalar product on $\lieg$ and use it to identify $\lieg^*$ with $\lieg$. Then we consider the \emph{norm square momentum map} $\parallel \mu \parallel^2:M\to \R$, where $\parallel\cdot\parallel$ is the norm induced by the $G$-invariant scalar product. 

Let $F$ be a connected component of the core and let $G_F:=\{g\in N(H):\, gF=F\}$. Then $F$ is a connected K\"ahler manifold, and the restriction $\mu_F:=\mu|_F$ is a momentum map for the $G_F$-action on $F$. In particular, the restriction of $\parallel \mu \parallel^2$ to $F$ is the norm square momentum map. The $G_F$-action on $F$ extends holomorphically to an action of $G_F^\C$. We point out that the good quotient of the set of the semistable points of the $G$-action on $M$ is homeomorphic to the good quotient of the set of the semistable points of the $G_F$-action on $F$. Under a mild assumption on the $G$-action, we prove that the set of semistable points of $M$ intersects $F$ in the set of the semistable points of $F$. If $M$ is compact, then $\parallel \mu \parallel^2$ induces a smooth stratification $\{S_\beta:\, \beta \in \mathcal B\}$ of $M$.  We prove that $\{S_\beta \cap F:\, \beta \in \mathcal B\}$ is the smooth stratification of $F$ induced by the restriction of $\parallel \mu \parallel^2$ to $F$.

\section{Preliminaries} \label{sect preliminaries}

In this section, we recall some basic definitions and facts about proper Lie group actions. We will also fix some notation that we shall use throughout the paper.\\

Let $G$ be a Lie group and $M$ be a smooth manifold. We say that $G$ \emph{acts} on $M$ \emph{on the left}, or that $M$ admits a \emph{left $G$-action}, if there exists a group homomorphism $A$ form $G$ to the group of diffeomorphisms of $M$, $\text{Diff}(M)$, such that the map
\[
G\times M \lra M, \qquad (g,x)\mapsto gx:=(A(g))(x),
\]
is smooth. The manifold $M$ is also called a \emph{$G$-manifold}. If, instead, $A$ determines a group anti-homomorphism, we call $A$ a \emph{right action}. Note that every right (left, resp.) action can be turned into a left (right, resp.) action by composition with the inversion map $G\to G$, $g\mapsto g^{-1}$. If not otherwise specified, we will always refer to left actions.

If $M$ and $M'$ are two $G$-manifolds, with $G$-actions $A$ on $M$ and $A'$ on $M'$, a map $f:M\to M'$ is called \emph{$G$-equivariant} if we have $f\circ A(g)=A'(g)\circ f$ for all $g\in G$. In particular, if $A'$ is the \emph{trivial action} on $M'$, i.e. $A'(g)=\text{id}_{M'}$ for all $g$ in $G$, a $G$-equivariant map $f:M\to M'$ is then said to be \emph{$G$-invariant}. 

Let $x\in M$. For a given $G$-action, the \emph{orbit through} $x$ is defined as the set $G(x):=\{gx:\, g\in G\}\subseteq M$ and the \emph{stabilizer} or the \emph{isotropy subgroup} of $x$ is the set $G_x:=\{g\in G\:\, gx=x\}\subseteq G$ which is clearly a closed subgroup of $G$, and so a Lie group \cite{DK}. If $G(x)=\{x\}$, equivalently $G_x=G$, then $x$ is called a \emph{fixed point} of the $G$-action on $M$.
We will often refer to more specific types of action. In particular, we will call a $G$-action on $M$

\begin{itemize}
    \item \emph{effective} if the kernel of $A$ is given by $\{e\}$, where $e$ is the identity element of $G$;
    \item \emph{free} if $G_x=\{e\}$ for all $x\in M$;
    \item \emph{transitive} if for any $x,y\in M$ there exists $g\in G$ with $gx=y$;
    \item \emph{proper} if the map $G\times M \lra M \times M,\  (g,x) \mapsto (gx,x)$, is proper in the topological sense, i.e., the preimage of any compact set is compact. 
\end{itemize}
\indent An important fact to recall is that if $G$ acts properly on $M$, then the stabilizer subgroups are compact and the $G$-orbits are embedded, closed submanifolds of $M$.  

Let $x,y\in M$. The condition $x\in G(y)$ defines an equivalence relation in $M$, that we denote by $x\sim y$. The quotient set $M/G:=M/\sim$, consisting of the collection of the $G$-orbits, is called the \emph{orbit space} of $M$.
Let $\pi:M\lra  M/G$ denote the natural projection. $M/G$ is naturally endowed with the quotient topology, which is the finest topology on $M/G$ that makes the projection $\pi:M\lra M/G$ continuous. It is easy to check that $\pi$ is also an open map.

 Throughout the paper we shall always assume that $M/G$ is connected. This means that $M$ could have many different connected components but the natural $G$-action on the set of the connected components of $M$ is transitive.\\


From now on, we will always consider proper actions. Two basic results for proper actions are the Slice Theorem and the Tubular Neighborhood Theorem (see \cite[Chapter 3]{BM}). We briefly review some of their consequences.

Let $x\in M$. Since $G_x$ is compact, there  exists a $G_x$-invariant Riemannian metric $\mathtt{l}$ on $M$. This metric, for instance, is defined by averaging an arbitrary Riemannian metric, $\langle \cdot,\cdot \rangle$, on $M$ using the $G_x$-action, that is,
\[
\mathtt l ( X,Y ) (p):=\int_{G_x} \langle (\mathrm d g)_p (X), (\mathrm d g)_p (Y)\rangle (gp) \, \tilde \omega,
\]
where $\tilde \omega$ is a right invariant volume form on $G_x$.
Let $(T_x G(x))^\perp$ denote the orthogonal complement of $T_x G(x)$ in $T_xM$ with respect to $\mathtt l$.  Then $$T_x M=T_x G(x) \oplus (T_x G(x))^\perp $$ is a  $G_x$-invariant orthogonal splitting. The Slice Theorem proves that there exists $\epsilon>0$ such that $\mathcal N_x=\exp_x (B_\epsilon \cap (T_x G(x))^\perp )$ is a submanifold of $M$ satisfying the following conditions:
\begin{itemize}
\item  $\mathcal N_x$ is $G_x$-invariant;
\item $G\mathcal N_x$ is an open subset of $M$;
\item if $g\in G$ and $y\in \mathcal N_x$ are such that $gy\in \mathcal N_x$, then $g\in G_x$.
\end{itemize}
A $G_x$-invariant submanifold of $M$ satisfying the above conditions is called a \textit{slice at} $x$ for the $G$-action on $M$.

The isotropy representation is the natural representation of $G_x$ on the tangent space $T_x M$ given by
\[
G_x \lra \mathrm{GL}(T_x M), \quad g \mapsto \mrm (d g)_x.
\]
If $M$ is connected, then the kernel of the isotropy representation coincides with the kernel of the homomorphism $A : G \lra \mathrm{Diff}(M)$. Indeed, it is well-known that if an isometry $f:M\lra M$ satisfies $f(x)=x$ and $(\mathrm d f)_x=\mathrm{id}_{T_xM}$ then $f=\mathrm{id}_M$; see \cite{do}. 

Since the isotropy representation preserves the orthogonal splitting $T_x M=T_x G(x) \oplus (T_x G(x))^\perp$, we can consider its restriction to $(T_x G(x))^\perp$, i.e. 
$G_x \lra \mathrm{O}((T_x G(x))^\perp ),\, g \mapsto {(\mrm{d} g_x)|}_{{(T_x G(x))^\perp}},$
which is called the \emph{slice representation}.

The compact Lie group $G_x$ acts on $G \times (T_x G(x))^\perp$ as follows
\[
h(g,v):=(gh^{-1},hv).
\]
This action is free and proper. Hence, its orbit space, denoted by  $E=\bigr(G\times (T_x G(x))^\perp\bigr) /G_x = G \times_{G_x} (T_x G(x))^\perp$,
is a homogeneous vector bundle associated to the principal bundle $G \lra G/G_x$; see \cite{BM, DK}. The group $G$ acts on $E$ as follows:
$h[g,v]:=[hg,v]$.  Since $\bigr(G ( [e,v])\bigr)\cap (T_x G(x))^\perp=G_x(v)$ one can check that the orbit space $E/G$ is homeomorphic to the orbit space $(T_x G(x))^\perp /G_x$. Let $\mathcal{N}_x$ be a slice at $x$. The Tubular Neighborhood Theorem states that the $G$-invariant neighborhood $\mrm{Tub}(G(x)):=G\mathcal{N}_x$ of $G(x)$ is $G$-equivariantly diffeomorphic to $G\times_{G_x}\mathcal{N}_x$ through the map 
\[
G\times_{G_x}\mathcal{N}_x\lra\mrm{Tub}(G(x)), \quad [g,v]\mapsto g\exp_x(v)
\]
 An important consequence of this is that the orbit space $M/G$ is locally compact, Hausdorff, second countable, and hence paracompact; see again \cite{BM, DK}.\\

Let now $\mathfrak g$ be the Lie algebra of $G$, let $x\in M$ be an arbitrary point, and let $\mathfrak{g}_x$ be the Lie algebra of the isotropy group $G_x$. Given $\xi\in\mathfrak{g}$, we define the \textit{infinitesimal generator of the action} or \emph{fundamental vector field} associated with $\xi$ as the map $\xi_M:M\lra TM$ given by 
 \begin{equation}\label{infgen}
 x\mapsto \xi_M(x)\coloneqq \tfrac d {dt}\bigr\vert_{t=0}\exp(t\xi)x\in T_xG(x),
 \end{equation}
 where $\exp:\mathfrak g\to G$ is the Lie exponential map of $G$. In fact,  one can show that $\xi_M$ is a smooth vector field on $M$, see \cite{BM}.
 We may identify $T_xG(x)$ with $\mathfrak g/\mathfrak g_x$. Indeed, the map $\mathfrak g \lra T_xG(x)$ given by  
\begin{equation*}
   \xi\mapsto \xi_M(x)\ , 
\end{equation*}
induces an isomorphism $\mathfrak g/\mathfrak g_x\cong T_xG(x)$. Under this identification, the $G_x$-action on $T_xG(x)$ corresponds to the adjoint action of $G_x$ on $\mathfrak g/\mathfrak g_x$. \\

Let $\mathcal U=\{U_\alpha \}_{\alpha \in I}$  be an open covering of $M$ consisting of $G$-invariant open subsets. The Tubular Neighborhood Theorem and, consequently, the paracompactness of $M/G$, prove there exists a $G$-invariant partition of unity
$\{ f_j :M \lra \R\}_{j\in J}$ subordinate to $\mathcal U$. This means each $f_j$ is a smooth $G$-invariant function, $f_j \geq 0$ and the support of $f_j$ is contained in some $U_{\alpha(j)} \in \mathcal U$. Moreover, $\{\mathrm{supp} f_j /G \}_{j\in J}$ forms a locally finite family of compact subsets of $M/G$, and $\sum_{j\in J} f_j =1$. Such partitions of unity can be used to glue together $G$-invariant structures. As an application (see \cite{BM, DK}), one can prove the existence of a $G$-invariant Riemannian metric $\mathtt{h}$ on $M$. The $G$-action is then called \emph{isometric} with respect to $\mathtt h$. We point out that there always exist a $G$-invariant Riemannian metric $\mathtt g$ in the conformal class of $\mathtt h$ which is complete. In fact, one has the following result by Podest\`a and Spiro \cite{ps}.
\begin{prop}\label{complete}
Let $G$ be a Lie group acting properly on $M$. If $M/G$ is connected, there exists a complete Riemannian metric $\mathtt g$ such that any element of $G$ is an isometry of $(M,\mathtt g)$.
\end{prop}
\begin{proof}
 We may assume, without loss of generality, that $M$ is connected. Indeed, let $M'$ be a connected component of $M$ and consider $G_{M'}:=\{g\in G:\, gM'=M'\}$. It is well known that the identity component of $G$, that we denote with $G^o$, is a closed normal subgroup of $G$ \cite{DK}. Since $G^o \subseteq G_{M'} \subseteq G$, it follows that $G_{M'}$ is a closed subgroup of $G$ and hence the $G_{M'}$-action on $M'$ is proper. Assume there exists a $G_{M'}$-invariant, complete Riemannian metric $\mathtt h'$ on $M'$. Let $M''$ be another connected component of $M$. Since $M/G$ is connected there exists $g_o \in G$ such that $g_o M''=M'$. Therefore, $G_{M''}=g_0^{-1}G_{M'}g_0$ and $\mathtt h'':=g_0^*\mathtt h'$ is a $G_{M''}$-invariant, complete Riemannian metric on $M''$. 


For any $x\in M$, we define
$r(x)$ to the supremum of the set of positive real numbers $r$ such that the open ball $B_r (x):=\{y\in M:\, d(x,y) < r\}$ is relatively compact. If $r(x)=+\infty$ for some $x\in M$, then any open ball is relatively compact. By the Hopf--Rinow Theorem  $(M,\mathtt h)$ is a complete Riemannian manifold \cite{do}. Hence, we may assume $r(x) <+\infty$ for any $x\in M$.

Let $x,y\in M$ be such that $y\in B_{r(x)} (x)$. Then $r(y)\geq r(x)-d(x,y)$ and so
\[
|r(x)-r(y)|\leq d(x,y).
\]
In \cite{no}, Nomizu and Ozeki proved the existence of a positive smooth real function $f:M \lra \R$ satisfying $f(x) \geq \frac{1}{r(x)}$ for any $x\in M$ and $(M,f^2 \mathtt h)$ is complete. Indeed, if we denote by  $d'$ the distance induced from  $\mathtt g$, then
If we denote by  $d'$ the distance induced from  $\mathtt g$, then
\[
\left\{ y\in M:\, d'(x,y) <\frac{1}{3} \right\} \subset \left\{y\in M:\, d(x,y) \leq \frac{r(x)}{2}\right\}.
\]
for any $x\in M$.
This implies that  any open ball of radius $\frac{1}{3}$ with respect to $\mathtt g$ is relative compact and so any Cauchy sequence with respect to $d'$ admits a convergent subsequence. By the Hopf--Rinow Theorem $(M,f^2 \mathtt h)$  is a complete Riemannian manifold \cite{do}.

Since $\mathtt h$ is $G$-invariant it follows that the function $r(\cdot)$ is a $G$-invariant continuous function. Our result follows form the existence of  a $G$-invariant function $f:M \lra \R$ such that $f(x) \geq \frac{1}{r(x)}$ for any $x\in M$. To show that such a function exists, we proceed as follows.

Let $x\in M$. By the Tubular Neighborhood Theorem, we may think a
$G$-invariant neighborhood of $G(x)$ as the homogeneous vector bundle $G\times_{G_x} \mathcal N_x$.
Let $B_n=\{v\in \mathcal N_x:\, \parallel v \parallel <n\}$.
The open covering $\mathcal U=\{GB_n\}_{n\in \N}$ of $G\times_{G_x} \mathcal N_x$ is $G$-invariant. Let $\{h_i:G\times_{G_x} \mathcal N_x \lra [0,1]\}_{i\in I}$ be a
$G$-invariant partition of unity subordinate to $\mathcal U$. The function
\[
\tilde f=\sum_{i\in I} \frac{h_{n(i)}}{r_{n(i)}},
\]
where $n(i)\in \N$ is such that $\mathrm{supp}\, h_i \subset GB_{i(n)}$ and $r_{n}=\mathrm{inf}_{z\in GB_n} r$,
satisfies $\tilde f \geq \frac{1}{r}$. Now, using the $G$-invariant partition of unity subordinated to $\{G\times_{G_x} \mathcal N_x \}_{x\in M}$,
one can construct a $G$-invariant smooth function $f$ satisfying
\[
f(y)\geq \frac{1}{r(y)},
\]
for any $y\in M$, concluding the proof.
\end{proof}
\begin{rem}
The above result is stated and proved in \cite[Teorema 3.4, page 101]{ps}. These lecture notes are available online, but have not been published yet. For this reason, we give a proof.
\end{rem}
From now on we always assume that $(M,\mathtt g)$ is a complete Riemannian manifold and that the $G$-action on $M$ is isometric with respect to $\mathtt g$, i.e. any element $g\in G$ is an isometry of $(M,\mathtt g)$. Let $d$ denote distance function defined by the Riemannian metric $\mathtt g$ \cite{do}. We remark that $d(x,y)=+\infty$ whenever $x$ and $y$ do not belong to the same connected component.

Let $\pi: M \lra M/G$ be the natural projection. For any given $x\in M$, we let $[x]:=\pi(x)\in M/G$. Given $p,q\in M$, we define 
\[
d_{M/G}([p],[q]):=\mrm{inf} \left\{d(gp,g'q):\, g,g'\in G \right\}.
\]
Note that $d_{M/G}([p],[q])< +\infty$ is due to the fact that $M/G$ is connected. As we will see below, $d_{M/G}$ is a distance function, called the \textit{orbital distance}.
Since the $G$-orbits are closed, the following result holds.
\begin{lemma}\label{distance-submetry}
$(M/G,d_{M/G})$ is a metric space and $\pi:M \lra M/G$ is submetry. This means that for any $p\in M$ and any
$\epsilon >0$, the $\epsilon$-ball $B_\epsilon (p)$ around $p$ is mapped onto $B_\epsilon ([p])$. In particular, $(M/G,d_{M/G})$ is a complete metric space.
\end{lemma}
\begin{proof}
Let $[p],[q]\in M/G$. Then $\dist{[p]}{[q]}\geq 0$. Assume $\dist{[p]}{[q]}=0$. Since 
\begin{equation}\label{d}
d_{M/G}([p],[q])=\mrm{inf} \left\{ d(p,gq):\, g\in G \right\},
\end{equation}
it follows that there existe a sequence $g_n \in G$ such that 
\[
d(p, g_n q ) \to 0.
\]
and so that $g_nq\to p$ in $M$. Since the $G$-action on $M$ is proper, it follows that $[p]=[q]$.

Let $[p],[q],[z]\in M/G$. For any $g,g'\in G$, we have
\[
\dist{[p]}{[q]}\leq d(p,gq)\leq d(p,g'z)+d(gq,g'z),
\]
and so
\[
\dist{[p]}{[q]}\leq d_{M/G} ([p],[z])+d_{M/G} ([z],[q])
\]
This implies the triangle inequality of the distance function on $M/G$.

Let $z\in B_\epsilon (p)$. Since $\dist{[p]}{[z]}\leq d(p,z)$, it follows that $\pi(B_\epsilon (p) )\subset B_\epsilon ([p])$.

Let $[q]\in  B_\epsilon ([p])$. By $(\ref{d})$ it follows there exists $g\in G$ such that $d(p,gq)<\epsilon$. This implies $\pi(B_\epsilon (p) )= B_\epsilon ([p])$.

Let $[p]\in M/G$ and let $\epsilon>0$. We claim that
$\overline{B_\epsilon ([p])}=\pi(\overline{B_\epsilon (p)})$. Therefore, the open balls in $(M/G,\mrm{d}_{M/G})$ are relatively compact.

Let $[z]\in \overline{B_\epsilon ([p])}$. Since $\pi$ is a submetry, there exists a sequence $\{y_n\}_{n\in \N}$ such that $y_n \in B_\epsilon (p)$ and
$
\pi(y_n)=[y_n] \to [z]$.  By the Hopf-Rinow Theorem (see \cite{do}), the open ball $B_\epsilon (p)$ is relatively compact. Hence we may assume that $y_n \to y_o \in \overline{B_p (\epsilon)}=\{z\in M:\, d(p,z)\leq \epsilon\}$. Therefore $\pi(y_o)=[z]$ and so  $\overline{B_\epsilon ([p])}\subseteq \pi(\overline{B_\epsilon(p)})$. The other inclusion is easy to check.

Let $\{ [p_n ] \}_{n\in \N}$ be a Cauchy sequence in $(M/G,d_{M/G})$ and let $\epsilon >0$. There exists $n_o$ such that $[p_n] \in B_\epsilon ([p_{n_o}])$ for any $n \geq n_o$.
Since $B_\epsilon([p_{n_o}])$ is relatively compact, it follows that $\{ [p_n ] \}_{n\in \N}$ has a convergent subsequence. This implies that the Cauchy sequence itself has a limit, concluding the proof.
\end{proof}

 One last fact we are going to recall is the following result by Kleiner, see \cite{BM} or Kleiner's PhD thesis \cite{Kl}. For the sake of completeness, we also present its proof. Let $p,q\in M$ and let $\gamma:[0,l]\to M$ be a geodesic with $\gamma(0)\in G(p)$ and $\gamma(l)\in G(q)$. We say that $\gamma$ is a \emph{Kleiner geodesic} from $G(p)$ to $G(q)$ if its length $L(\gamma)$ equals the distance $d\left(G(p),G(q)\right)$; it is obvious that a Kleiner geodesic is minimizing. Kleiner's result tells us the following.
\begin{lemma}\label{geodistanza}
Given any two points $p,q\in M$, there exists a Kleiner geodesic $\gamma:[0,l]\to M$ from $G(p)$ to $G(q)$. Moreover, for any such Kleiner geodesic we have $G_{\gamma(t)}=G_\gamma$ for any $t\in (0,l)$, where
\[
G_\gamma\coloneqq \{g\in G\ :\  g\gamma(s)=s \quad \text{for all } \;s\in[0,l]\}
\]
\end{lemma}
\begin{proof}
Let $x_n \in G(p)$ and $y_n \in G(q)$ be such that
\[
\lim_{n\to +\infty} d(x_n,y_n)=d(G(p),G(q))=l.
\]
By $(\ref{d})$ we may assume that $x_n=p$ and so $y_n $ lies in the same connected component of $p$. By the Hopf-Rinow Theorem there exists a minimizing geodesic $\gamma_n $, parameterized with the arc length, joining  $p$ and $y_n$.
The geodesic $\gamma_n :[0,t_n] \lra M$ is given by $\gamma_n (t)=\exp_p (t\xi_n)$, where $\parallel \xi_n  \parallel =1$ and $\gamma(t_n)=y_n$.  In particular, $L(\gamma_n)=t_n$ and
$\lim_{n\to +\infty } t_n = d(G(p),G(q))=l$. Up to passing to a subsequence, we may assume that the sequence $\{ \xi_n \}_{n\in \N}$ converges to some unit vector $v \in T_p M$. Now, keeping in mind that $G(q)$ is closed, it is easy to check that the geodesic $\gamma (t):=\exp_p(tv)$ is a minimizing geodesic joining
$p$ and $\gamma(r)\in G(q)$ and satisfying $L(\gamma)=d(G (p),G (q))$.

Let $t\in (0,r)$ and let $g\in G_{\gamma(t)}$. The curve
\[
\tilde \gamma (s)=\left\{\begin{array}{lc} g \gamma (s) & 0\leq s\leq t \\ \gamma(s) & t<s\leq r \end{array}\right.
\]
joins $g (\gamma (0))$  and $\gamma (r)$. Since $L(\tilde \gamma )=L(\gamma)=d(G(p),G(q))$ it follows that $\tilde \gamma$ is a minimizing geodesic, hence smooth. This implies $(\mrm{d} g)_{\gamma (t)}(\fis{\gamma}(t))=\fis{\gamma }(t)$. Since $g(\gamma(t))=\gamma(t)$, it follows that $g\in G_\gamma\coloneqq \bigcap_{s\in [0,r]}G_{\gamma(s)}$, concluding the proof.
\end{proof}

\section{Reduction to the core} \label{sect reduction to core} 
In this section we consider a proper action of a Lie group $G$ on a smooth manifold $M$ and describe the construction of the core $\core$, following \cite{GS} and \cite{SS}. We assume that $M/G$ is connected, that $(M,\mathtt g)$ is a complete Riemannian manifold, and that $G$ acts on $(M,\mathtt g)$ by isometries. \\

Let $K$ be a compact subgroup of $G$ and let
\[
M^{K}:=\{p\in M:\, kp=p,\, \text{ for all } k\in K\},
\]
be the fixed point set of $K$.
By a Theorem of Koszul, the connected components of $M^{K}$  are closed embedded submanifolds of $M$, of possibly different dimensions, whose tangent space at $x\in M^K$ is given by
 $$T_x M^K=(T_x M)^K:=\{v\in T_x M:\, (\mrm d k)_x (v)=v \text{ for all }k\in K\},$$
see for instance \cite[Proposition 3.93]{BM}. Therefore,  if $Y$ is a connected component of $M^K$, then $Y$ is totally geodesic and $(Y,{\mathtt g|}_{Y})$ is a complete Riemannian manifold. 

 Let $N(K)$ be the normalizer of $K$ in $G$; i.e. $N(K)\coloneqq \{g\in G\ :\  gKg^ {-1}=K\}$. We have that
 $N(K)$ is a closed subgroup of $G$, and so a Lie group \cite{DK}. $N(K)$ preserves $M^K$ with $K$ as ineffective kernel. In particular, the $N(K)$-action on $M^K$ is proper.

 We now define $$M_{K}:=\{p\in M:\, G_p =K\}\ ,$$ and $$M_{(K)}:=\{p\in M:\, (G_p)=(K)\}\ ,$$ where $(K)$ the conjugacy class of $K$ in $G$. Since $G_{gp}=gG_pg^{-1}$, $M_{(K)}$ is the union of the all $G$-orbits of orbit type $(K)$. 
 Note that if $gM_{K} \cap M_{K} \neq \emptyset$ then $g\in N(K)$ and
$M_{K}$ is $N(K)$-invariant. Since $M_{(K)}=GM_{K}$, it is $G$-invariant. The following result is well known \cite{DK}, however, since it will be of use in the following, we recall its proof.
\begin{prop}\label{otl}
The connected components of $M_{K}$ and $M_{(K)}$ are embedded submanifolds of possibly different dimensions.
\end{prop}
\begin{proof}
Let $q\in M_{K}$. Let $E=G \times_{K} \mathcal N_q$ be the tube. 
Then $[g,v]\in M_{K}$ if and only if $gK_{v}g^{-1}=K$ and so if and only if $K_{v}=K$ and $g\in N(K)$. Therefore,
\begin{align}\label{tube in M_H}
M_{K} \cap E= N(K) \times_{K} \mathcal N_q^{K},
\end{align}
and so
\begin{align}\label{tube in M_(H)}
M_{(K)} \cap E= G \times_{K} \mathcal N_q^{K}.    
\end{align}
This concludes the proof.
\end{proof}
\begin{prop}\label{red1}
Both the $N(K)$-action on $M_K$ and the $G$-action on $M_{(K)}$ are proper and the natural inclusion $M_{K} \hookrightarrow M_{(K)}$ induces a diffeomorphism from $M_{K} / N(K)$ onto $M_{(K)} /K$.
\end{prop}
\begin{proof}
Since $M_{(K)}$ is $G$-invariant, it follows the $G$-action on $M_{(K)}$ is proper. Since $N(K)$ is closed in $G$, the $N(K)$-action on $M_{K}$ is proper as well.
Let $p\in M_{K}$. It is easy to check that  
\begin{align}\label{N(H)-princ.orbits}
	G(p)\cap M_K=N(K)(p).
\end{align}
Hence the inclusion $i:M_{K} \hookrightarrow M_{(K)}$ induces a continuous bijective map $\tilde i: M_{K} / N(K) \lra M_{(K)} /G$. By (\ref{tube in M_H}) and (\ref{tube in M_(H)}) it follows that $\tilde i$ is a diffeomorphism.
\end{proof}

\begin{prop}\label{fibration}
The map $S:G \times_{N(K)}  M_K \lra M_{(K)}$ given by $[g,q] \mapsto g q$ is a $G$-equivariant diffeomorphism.
\end{prop}
\begin{proof}
By (\ref{tube in M_H}) and (\ref{tube in M_(H)}) the map $S$ is a local diffeomorphism. Since $M_{(K)}=GM_{K}$ it follows that $S$ is surjective. Now we show that $S$ is injective. If $q,w\in M_K$ and $g,h\in G$ with $gq=hw$, then $h^{-1}g \in N(K)$ and so $[g,q]= [g (g^{-1} h ) , h^{-1} gq]=[h,w]$. This proves that $S$ is a diffeomorphism.
\end{proof}

 Since $M/G$ is connected, the Principal Orbit Theorem \cite{BM} states that there is a unique isotropy type $(H)$ such that the set of points of $M$ with orbit type $(H)$ is open and dense. Let $H$ be a principal isotropy and let
\[
\princ:=\{p\in M\ :\  (G_p)=(H)\},
\]
the $G$-invariant  open and dense subset of $M$ of all the $G$-orbits of orbit type $(H)$.  If $x\in \princ$, then $\dim G(x) \geq \dim G (y)$ for any $y\in M$ and $G(z)=G(x)$ for all $z\in\princ$. Furthermore, the slice representation at $x$ is trivial. By definition, the \emph{cohomogeneity} of the $G$-action on $M$ is the codimension of a principal orbit and it is denoted by $\mrm{cohom}(G,M)$. Moreover, $\princ /G$ is a connected manifold and the restriction  $\pi: \princ \lra \princ/G$ is a Riemannian submersion and a fiber bundle with typical fiber $G/H$. We define the dimension of $M/G$, which we denote by $\dim (M/G)$, as the dimension of $\princ /G$, which, in turn, coincides with the cohomogeneity of the $G$-action on $M$ \cite{BM}. 

Let $H$ be a principal isotropy and let $N\coloneqq N(H)$ denote the normalizer of $H$ in $G$. Since the slice representation at any $x\in M_H$ is trivial, it follows from the Tubular Neighborhood Theorem that $M_H$ is an $N$-invariant open subset of $M^H$ and $\princ = GM_H$. Keeping in mind that $H$ is a principal isotropy, Proposition \ref{fibration} shows that $\princ$ is a homogenenous $M_H$-bundle over $G/N$. Together with Proposition \ref{red1}, we get:

\begin{prop}\label{Reduction on principal stratum}
	The function
	\begin{align*}
		S&:G\times_{N} M_H \lra \princ, & [g,p]\mapsto gp, 
	\end{align*}
	is a $G$-equivariant diffeomorphism and the inclusion $i:M_H\hookrightarrow \princ$ induces a diffeomorphism $M_H/N\lra \princ /G$.
\end{prop}
	
 Since $M/G$ is connected, by the Principal Orbit Theorem $\princ/G$ is connected and thus by Proposition \ref{Reduction on principal stratum} also $M_H/N$ is connected. In particular, the connected components of $M_H$ have the same dimension. As a consequence of Lemma \ref{geodistanza}, we get the following proposition. 
 
\begin{prop}\label{red-isometry-principal}
	The inclusion $M_H\hookrightarrow \princ$ induces an isometry $M_H/N\lra \princ /G$ with respect to the orbital distances.
\end{prop}
\begin{proof}
 	Let $x,y\in M_H$. We claim that 
 	\begin{align} \label{distanceEq. principal}
 		d\bigr(N(x),N(y)\bigr)\coloneqq \mrm{inf} \{d(x,hy)\ :\  h\in N\} = d\bigr(G(x),G(y)\bigr),
 	\end{align}
 	the result will follow. Indeed, since $N$ is a subgroup of $G$, we have that $	d\bigr(N(x),N(y)\bigr)\geq d\bigr(G(x),G(y)\bigr)$. By Lemma \ref{geodistanza}, there exists a Kleiner geodesic $\gamma:[0,l]\lra M$ from $G(x)$ to $G(y)$ such that $\gamma(0)=x$. Note that $G_x\supseteq G_\gamma \subseteq G_{\gamma(l)}.$ Since $G_x=H$ is principal, we have that $G_\gamma=H$. Keeping in mind that $(G_y)=(G_{\gamma(l)})=(H)$, we obtain that $G_{\gamma(l)}=H$ and thus $\gamma(l)\in G(y)\cap M_H$. By (\ref{N(H)-princ.orbits}) we have that $\gamma(l)\in N(y)$. This shows that (\ref{distanceEq. principal}) holds and completes the proof.
\end{proof}

Let $C$ be a connected component of $M_H$ and let $G_C=\{g\in G \ :\  gC=C\}$. One can readily check that $N^0\subseteq G_C\subseteq N$ and so $G_C$ is a closed subgroup of $N$ given by a disjoint union of connected components of $N$. In particular, the $G_C$-action on $C$ is proper. Since $M_H/N$ is connected, and $\princ =GM_H$, we have that $\princ=GC$. Since the $G$-action on $M$ is by isometries and $C/G_C=M_H/N$, we have the following result.

\begin{prop}
	The function
	\begin{align*}
		S_C&:G\times_{G_C} C \lra \princ, & [g,p]\mapsto gp, 
	\end{align*}
	is a $G$-equivariant diffeomorphism and the inclusion $i:C\hookrightarrow \princ$ induces an isometry $C/G_C\lra \princ /G$ with respect to the orbital distance.
\end{prop}

Let $\prescript{}{c}{M}$ denote the closure of $M_H$ in $M$. We will refer to $\core$ as the \textit{core} of the $G$-action on $M$. Clearly, we have that
\[
M_H\subseteq \core\subseteq M^H
\]
and all three of these subsets are $N$-invariant. 

We now study some of the properties of $\core$. We begin with the following lemma.
\begin{lemma}\label{gvh}
    Let $F$ be a connected component of $M^H$ that meets $M_H$. Then for all $x\in F$, there exists $z\in \mathcal N_x$ such that $G_z=H$.
\end{lemma}

\begin{proof}
    Let $G_F\coloneqq \{g\in G\ :\  gF=F\}$. Clearly, $G_F$ is a non-empty open and closed subgroup of $N$. Since $N$ is a closed subgroup of $G$, the $G_F$-action on $F$ is proper. Let $x\in F$, since $F\cap M_H\neq \emptyset$, by Lemma \ref{geodistanza}, there is a Kleiner geodesic $\gamma:[0,l]\lra F$ from $G_F(x)$ to $G_F(y)$ with $\gamma(0)=x$ and $\gamma(l)=y\in F\cap M_H$. Since $G_y=H$ is principal, we have that $(G_F)_\gamma=(G_F)_y=H$. Now, by (\ref{tube in M_H}), $\gamma$ meets $G(y)$ orthogonally at $\gamma(l)=y$. It is well known that if a geodesic meets a $G$-orbit orthogonally, then it is orthogonal to every orbit it meets. Hence $\gamma$ meets $G(x)$ orthogonally at $x$ and therefore $v=\gamma'(0)\in (T_xG(x))^\perp$. It follows that for $\tau>0$ small enoguh, $z=\exp_x(\tau v)\in \mathcal N_x$. Hence, $G_z\subset G_x$. Thus, every element $g\in G$ fixing $z$ must also fix $x$ and, by the properties of the Riemannian exponential, must also fix $\gamma(l)$. We conclude that $G_z=H$. 
\end{proof}

\begin{rem}\label{core i union of conn.comps}
    Thus, the core $\core$ coincides with the set of points $x\in M$ such that $H$ is a principal isotropy of the slice representation $G_x\lra \mathrm{O}\bigr(\bigr(T_xG(x))^ \perp\bigr)$; equivalently,
	\begin{align*}
		\core = \{x\in M \ :\  \text{ there exists } z\in \mathcal N_x \text{ such that } G_z=H\}.	
	\end{align*}

    It also follows that $\core$ is the union of the connected components of $M^H$ that meet $M_H$. In particular, $\core$ is a closed embedded totally geodesic submanifold of $M$ and $\dim (\core) = \dim (M/G) + \dim (N/H)$. Hence $(\core,\mathtt g)$ is a complete Riemannian manifold and $(\core/N,d_{\core/N})$ is a complete metric space such that $\pi:\core\to \core/N$ is a submetry. We also have that

\label{slice-reduction}
        \begin{itemize}
            \item $M_H$ is open and dense in $\core$;
            \item The group $N$ acts transitively on the connected components of $\core$. In particular, $\core/N$ is connected;
            \item If $g\in G$ is such that $g(\core)=\core$, then $g\in N$.
        \end{itemize}

    Another easy observation is that if we let $x\in \core=\overline{M_H}$, then, by the identification $(\mathfrak g/\mathfrak g_x)\simeq T_xG(x)$, we have $\bigr(\mathfrak g/\mathfrak g_x)^H= \mathfrak n/\mathfrak g_x$, where $\mathfrak n$ is the Lie algebra of the normalizer $N=N(H)$. Hence $\bigr(T_xG(x)\bigr)^H=\bigr(T_xN(x)\bigr)$ and therefore,
 
 \begin{itemize}
     \item  The slice for the $N$-action on $\core$ at $x$ is given by $(\mathcal N_x)^H$, the fixed point set for the $H$-action on $\mathcal N_x$. 
 \end{itemize}
\end{rem}

Now we state and prove an extension of Proposition \ref{red-isometry-principal} to the core $\core$. If $G$ is compact, this result is already pointed out in \cite{BH} and used in \cite{Straume2}. The proof is given in the unpublished paper \cite{SS}, see also \cite{GS}.

\begin{thm} \label{reduction-global}
    The inclusion $\core\hookrightarrow M$ induces an isometry $\Theta:\core/N\lra M/G$.
\end{thm}

\begin{proof}
    As mentioned above, $M_H/N$ is open and dense in $\core/N$. Simlarly, $M_{\mrm{princ}}/G$ is open and dense in $M/G$. By Proposition \ref{red-isometry-principal}, we have that the induced map $\Theta:\core/N\to M/G$ restricts to an isometry $\Theta\bigr\vert_{M_H/N}:M_H/N\to M_{\mrm{princ}/G}$. Since both $(M/G,d_{M/G})$ and $(\core/N,d_{\core/N})$ are complete metric spaces, it is easy to check that $\Theta$ is the unique isometry extending $\Theta\bigr\vert_{M_H/N}$.
\end{proof}

\begin{rem}
	It follows easily from Theorem \ref{reduction-global} that for all $x\in \core$ we have 
	\[
		G(x)\cap \core = N(x).
	\]
	Let us give an alternative proof of this fact. Let $g\in G$ and suppose that $gx\in \core$. By Lemma \ref{gvh}, there is some $z\in \mathcal N_x$ such that $G_z=H$. Since $g\mathcal N_x=\mathcal N_{gx}$, we have that $w=gz\in \mathcal N_{gx}$ is a principal point for the $G_{gx}$-action on $\mathcal N_{gx}$. Again, by Lemma \ref{gvh}, there exists some $\theta\in G_{gx}$ such that $G_{\theta w}=H$. Now, let $h=\theta g$, hence $hx=gx$. We claim that $h\in N$. Indeed,
	\begin{align*}
		hHh^{-1}=\theta gG_{z}g^ {-1}\theta^{-1}=\theta G_w \theta^{-1}=G_{\theta w}=H.
	\end{align*}
\end{rem}

Theorem \ref{reduction-global} also has the following useful consequence.

\begin{lemma}\label{lemma fixed point rep}
    Let $V$ be a finite dimensional representation of a compact Lie group $G$. Let $H$ be a principal isotropy. The fixed point-set of $N=N(H)$ coincides with the fixed point set of $G$; i.e. $V^{N}=V^G$.
\end{lemma}
\begin{proof}
    Since $N$ is a subgroup of $G$, we have that $V^G\subseteq V^N$. Let us now prove the reverse inclusion. Fix an inner product on $V$ for which the $G$-action is orthogonal. Consider the $G$-action on the unit sphere $M=\mathbb S^l$ of $V$. Clearly, the core $\core$ is given by the unit sphere $\mathbb S^r$ of $V^H$. Note also that if $x\in \core$ is a fixed point for the $N$-action, then so is $y=-x\in\mathbb S^r=\core$. By Theorem \ref{reduction-global}, we have that the natural map $ \core/N\lra M/G$ is an isometry and hence
    \[
    \pi = d_{\core}(x,y)=d_{\core/N}\bigr(N(x),N(y)\bigr)=d_{M/G}\bigr(G(x),G(y)\bigr)=\inf_{g\in G}\;d_{M}(gx,y).
    \]
    Note that in the sphere, any two non-antipodal points have distance $<\pi$. It follows that $x$ is a fixed point of the $G$-action on $M$ and thus a fixed point of the $G$-action on $V$.
\end{proof} 

Let us now study some properties of the stratification of $\core/N$. We begin with the following remark.

For any compact subgroup $K\subseteq G$, the connected components of the set $M_{(K)}$ of points of orbit type $(K)$ are embedded submanifolds. The partition of $M$ into these sumbanifolds provides a Whitney Stratification of $M$ (See e.g. \cite[Chapter 2.7]{DK}). The partition of $M/G$ by the connected components of the smooth manifolds $M_{(K)}/G$ also gives a Whitney stratification of $M/G$ called the \emph{orbit type stratification of }$M/G$; each stratum is called an \emph{orbit type stratum}. Let $x\in M$, and let $K=G_x$. The Tubular Neighborhood Theorem implies that $\mathcal N_x/K$ gives an open neighborhood of $[x]\in M/G$, where $[x]$ denotes the image of $x$ under the projection $\pi:M\lra M/G$. For $y\in \mathcal N_x$, $[y]$ belongs to the same orbit type stratum of $x$ if and only if $y\in (\mathcal N_x)^K$; i.e. if $y$ is a fixed point for the $K$-action on $\mathcal N_x$.

Now assume $x\in \core$, and recall that the $K=G_x$-action on $\mathcal N_x$ is equivalent to the slice representation $K\lra \mathrm O(V)$, where $V=\bigr(T_xG(x)\bigr)^\perp$. As discussed above, the slice representation for the $N$-action on $\core$ at $x$ is the subrepresentation $N_x=N_K(H)\lra \mathrm O(V^H)$. Let $y\in \mathcal N_x$, Lemma \ref{lemma fixed point rep} now implies the following:
\begin{itemize}
    \item The image of $y$ in $\core/N$ belongs to the same $N$-orbit type stratum of $[x]$ if and only if its image in $M/G$ belongs the same $G$-orbit type stratum of $[x]$.
\end{itemize}
In summary, we have proved the main theorem of this section.

\corereduction*

\begin{cor} \label{ChevCore}
    Let $G$ act properly on $M$ with $M/G$ connected; let $H$ be a principal isotropy. The restriction $C^\infty(M)\to C^\infty(\core)$ induces an isomorphism of algebras 
    \begin{align*}
        C^\infty(M)^G\to C^\infty(\core)^{L}.
    \end{align*}    
\end{cor}
\begin{proof}
    Recall that a continuous map $f:M/G\lra \mathbb R$ is said to be smooth if and only if $f$ restricts to a smooth map on every stratum of $M/G$. It is well-known that the pullback $\pi^*$ of the quotient map $\pi:M\lra M/G$ defines an isomorphism $C^\infty(M/G)\lra C^\infty (M)^G$ (see e.g. \cite[Chapter 2]{DK} ). By Theorem \ref{stratifiediso}, the inclusion $\iota:\core \hookrightarrow M$ induces an isomorphism $\core/L\lra M/G$ of stratified spaces. Therefore, the restriction $\iota^*:C^\infty(M)\lra C^\infty (\core)$ induces an isomorphism $C^\infty(M)^G\lra C^\infty(\core)^L$.
\end{proof}
Let $F$ be a connected component of the core $\core$ and let $G_F:=\{g\in G:\, gF=F\}$. We have that  $N^o \subseteq  G_F\subseteq N$ and $M_H \cap F$ is the set of principal points for the $G_F$-action on $F$. Clearly, for any $p\in F$ we have 
\[
G(p) \cap F=G_F (p). 
\]
We have proved $G(\core)=M$. Since $N$ acts transitively on the set of the connected components of the core, it follows that $GF=M$. Therefore, the inclusion $F \hookrightarrow M$ induces a bijective continuous map 
\[
\Theta_F : F/G_F \lra M/G.
\]
Let $p,q\in F$. By Lemma \ref{geodistanza} there exists $g\in N$ and a Kleiner geodesic $\gamma:[0,l] \lra \core$ joining $p$ and $gq$, contained in $F$, such that 
\[
d(N(p), N(q))=L(\gamma).
\]
This means $gq\in F$ and so $g\in G_F$. In particular, 
\[
d(N (p), N (q))=d(G_F (p),G_F (q)).
\]
We have thus proved the following theorem.
\begin{thm}\label{reduction2}
Let $F$ be a connected component of $\core$. Then the inclusion $F \hookrightarrow M$ induces an isometry $\Theta_F : F/G_F \lra M/G$.
\end{thm}

\section{Reduction principle for polar actions} \label{sect polar}
Let $G$ be Lie group acting properly and isometrically on a complete Riemannian manifold $(M,\mathtt h)$ with $M/G$  connected. The $G$-action on $M$ is called \emph{polar} if there exists a complete, connected, immersed submanifold $\Sigma$ of $M$ which intersects all the orbits orthogonally \cite{Path1,Path}. Such submanifold $\Sigma$ is called a \emph{section}.  If there is a section which is flat with respect to the induced Riemannian metric, then the action is called \emph{hyperpolar}.

Let $\overline{M}$ be the connected component of $M$ containing $\Sigma$. Since $G$ acts transitively on the set of the connected components of $M$, the following remark holds.
\begin{lemma}
The $G$-action on $M$ is polar if and only if the $G^o$-action on $\overline{M}$ is polar and so if and only if the $G^o$-action on any connected component of $M$ is polar.
\end{lemma}
We shall prove that a reduction principle holds for polar actions; this result is known by the experts. For the sake of completeness we give a proof. 

We shall use the following lemma, the proof of which can be found in \cite{Gorodski}.
\begin{lemma}
Let $x\in M$ and let $\mathcal V_x =\exp_x ((T_x G(x))^\perp )$.   Then
\begin{itemize}
\item $\mathcal V_x $ meets all the orbits of $G$;
\item if the $G$-action is polar and $x\in \princ$, then $\Sigma=\mathcal V_x$ is the unique section  through the point $x$.
\end{itemize}
\end{lemma}
Let $H$ be a principal isotropy and let $L=N/H$ and consider the $L$-action on the core $\core=\overline{M_H}$.



\polarred*
\begin{proof}
Since the actions of $L$ and $N$ on $\core$ have the same orbits, it suffices to prove the theorem for either the $N$-action or the $L$-action on $\core$.

Assume that the $G$-action on $M$ is polar. Let $x\in (\core)_{\mrm{princ}}$. Since $x\in \princ$, there exists a unique section through $x$ for the $G$-action on $M$ and it is given by $\Sigma=\mathcal V_x$. Therefore, keeping in mind that $H$ acts trivially on $\mathcal V_x$,  it follows that $\Sigma \subseteq \core$. Since $N\subset G$ and $L$-orbits coincide with $N$-orbits, it follows that $\bigr(T_yG (y)\bigr)^ \perp\subset \bigr(T_yL( y)\bigr)^\perp$ for all $y\in \core$. This shows that $\Sigma$ is a totally geodesic immersed submanifold of $\core$ meeting every $L$-orbit orthogonally. Thus $L$ acts polarly on $\core$.

Conversely, assume that the $L$-action on $\core$ is polar. Let $x\in (\core)_{\mrm{princ}}$. Let $\Sigma$ denote the unique section through $x$ for the $L$-action on $\core$. Recall that $(\core)_{\mrm{princ}}=M_H$ and hence by (\ref{tube in M_H}) there is an $N$-invariant neighborhood of $x$ in $\core$ that is $N$-equivariantly isometric to $N/H\;\times  \mathcal N_x$. It follows that $\Sigma=\mathcal V_x$. Since $\core$ is a totally geodesic submanifold of M, it follows that $\Sigma$ is a totally geodesic submanifold of $M$. Note that by Theorem \ref{reduction-global} we have that $\Sigma$ meets every $G$-orbit. Let $y\in \Sigma$. Keeping in mind that the slice at $y$ for the $N$-action on $\core$ is contained in the slice at $y$ for the $G$-action on $M$, it follows that $T_y\Sigma$ is orthogonal to $T_yG(y)$.  Since the $\core$ is totally geodesic, the reduction principle holds for hyperpolar actions as well. This completes the proof.   
 \end{proof}
\begin{cor}
Let $F$ be a connected component of $\core$. Then $N^o$ acts polarly, respectively hyperpolarly,  on $F$ if and only if $G$ acts polarly, respectively hyperpolarly, on $M$.
\end{cor}

Note that in the proof of Theorem \ref{reduction-polar} we have that if the $G$-action on $M$ is polar then there is a section $\Sigma$ contained in $\prescript{} c M$. Palais and Terng \cite{Path1} studied polar actions with the property that $N/H$ is discrete. In this case we have that any connected component of the core is a section for the $G$-action on $M$.
\begin{defini}
The $G$-action on $M$ is called infinitesimally polar if all the slice representations are polar.
\end{defini}
It has been shown by Lytchak and Thorbergsson in \cite{LT}
that the orbit space of a proper and isometric action of a Lie group on a Riemannian manifold,
with the quotient metric space structure, is a Riemannian orbifold if and only if the action is infinitesimally polar. As a consequence, the following reduction principle holds.

\begin{cor}\label{infpolar-reduction}
The $G$-action on $M$ is infinitesimally polar if and only the $N$-action on $\core$ is infinitesimally polar.
\end{cor}
\begin{proof}
Since $M/G$ is isometric to $\core/N$ the result follows directly by the Theorem of Lytchak and Thorbergosson. The other proof involves the result proved in this paper. 

Let $x\in \core$. The slice representation of $G_x\cap N$ on $T_x N(x)^\perp\cap T_x\left(\core\right)$ is the core of the slice representation of $G_x$ on $T_x G(x)^\perp$. By Theorem \ref{reduction-polar}, the $G_x\cap N$-action on $T_x N(x)^\perp\cap T_x\left(\core\right)$ is polar if and only if the $G_x$-action on $T_x G(x)^\perp$ is polar. Finally, one can easily verify that if the slice representation at $x$ is polar, so is the slice representation at $gx$ for any $g\in G$. Since $\core$ intersects all $G$-orbits of $M$, this finishes the proof.
\end{proof}
\begin{defini}
A $k$-section for the $G$-action on $M$, where $k$ is a nonnegative integer,  is a connected, complete submanifold $\Sigma$ of $M$ such that the following hold:
\begin{enumerate}
\item $\Sigma$ is totally geodesic;
\item  $\Sigma$ intersects all the $G$-orbits;
\item for every principal point $p\in \Sigma$, we have that $T_p \Sigma$ contains the normal space $(T_p G( p) )^\perp$ as a subspace of codimension $k$;
\item  for every principal point $p\in \Sigma$ we have that if $gp\in \Sigma$ for some $g\in G$, then $g\Sigma=\Sigma$.
\end{enumerate}
The integer
\[
\mrm{copol}(G,M):=\mathrm{min}\{k\in \mathbb N:\, \mrm{there\ exists\ a\ }k\text{-}\mrm{section}\ \Sigma \subset M\},
\]
is called the \emph{copolarity} of the $G$-action on $M$. If $\Sigma$ is a $\mrm{copol}(G,M)$-section then we say that $\Sigma$ is minimal.
\end{defini}
The notion of copolarity was introduced in \cite{got}. It is obvious that $\mrm{copol}(M,G)=0$ corresponds to a polar action. 

\begin{rem}\label{core_is_a_k_section}
    Note that the core $\core$ gives a $k$-section with $k=\dim (N/H)$. Therefore, we have that $0\leq \mrm{copol}(G,M)\leq \dim (N/H)$.
\end{rem}
\begin{thm}\label{copolarity-reduction}
$\mrm{copol}(M,G)=\mrm{copol}(L,\core)$.
\end{thm}
\begin{proof}
Let $\Sigma$ be a minimal $k$-section for the $G$-action on $M$ passing through $x\in (\core)_{\mrm{princ}}$. By Remark \ref{core_is_a_k_section}, we have that $\Sigma\subseteq \core$. Since $\Sigma$ meets all the $G$-orbit, it also meets all the $L$-orbits.  Keeping in mind that  $(T_x L(x))^\perp\cap T_x(\core)=(T_x G(x))^\perp$, it follows that (3) is satisfied. Property (4) also holds and therefore $\Sigma$ is a $k$-section for the $L$-action on $\core$.

Conversely, let $\Sigma$ be a $k$-section for the $L$-action on $\core$. By Theorem \ref{reduction2}, $\Sigma$ meets all the $G$-orbits. Since both item (1) and (3) follow directly by the fact that $\Sigma$ is a $k$-section for the $L$-action on $\core$, it suffices to prove item (4).

Let $p\in \Sigma $ be a principal point such that $gp\in \Sigma$ for some $g\in G$. Since $G (p)\cap \core =L(p)$, it follows that there exists $n\in N$ such that $g=nh$ for some $h\in H$. Therefore, item (4) holds, concluding the proof.
\end{proof}

\section{Coisotropic Actions}\label{sa}
Let $(M,\omega)$ be a symplectic manifold, with symplectic form $\omega$. In the following, we do not assume that $M$ is connected but possibly
has many different connected components of the same dimension.

We start by recalling some basic definitions and facts. 

Let $P$ be a submanifold of $M$ and let $x\in N$. The \textit{symplectic orthogonal} of $T_xP\subseteq T_xM$ is defined as
\[
T_x P^\omega:=\{v\in T_x M:\, \omega_x(v,w)=0,\, \forall w\in T_x N\}.
\]
\begin{defini}
Let $P$ be a submanifold of $(M,\omega)$. We say that
\begin{enumerate}
\item $P$ is an \textit{isotropic submanifold} if at any point $x\in P$, $T_x P$ is an \textit{isotropic subspace} of $T_x M$, i.e., 
$T_x P\subseteq T_x P^\omega$;
\item $P$ is a \textit{coisotropic submanifold} if at any point $x\in P$, $T_x P$ is a \textit{coisotropic subspace} of $T_x M$, i.e.,
$T_x P^\omega \subseteq T_x P$;
\item $P$ is a \textit{Lagrangian submanifold} if at any point $x\in P$, $T_x P$ is a \textit{Lagrangian subspace} of $T_x P$, i.e.,
$T_x P=T_x P^\omega$;
\item $P$ is a \textit{symplectic submanifold} if at any point $x\in N$, $T_x P$ is a \textit{symplectic subspace} of $T_x M$, i.e., $T_x P\cap T_x P^\omega =\{0\}$;
\end{enumerate}
\end{defini}

Recall that a \textit{symplectomorphism} of $(M,\omega)$ is a diffeomorphism $\phi:M \lra M$ such that $\phi^* \omega=\omega$. 
A vector field on $M$ is called a \emph{locally Hamiltonian field} if the Lie derivative of $\omega$ along $X$ is zero: $\lie X \omega =0$. By Cartan's formula \cite{Ca}, we have
\[
\lie X \omega=\mathrm{d} \it{i}_X \omega + \it{i}_X \mathrm d \omega= \mathrm{d} \it{i}_X \omega,
\]
and so $X$ is a locally Hamiltonian vector field of $M$ if and only if $\contr_X \omega$ is closed. In particular, if $\contr_X \omega =\text{d} f$, for some $f\in C^\infty(M)$, we call $X$ the \textit{Hamiltonian vector field} associated with $f$, or the \textit{symplectic gradient} of $f$.
Let $f,g\in C^\infty (M)$. The \textit{Poisson bracket} of $f$ and $h$ is defined as
$\{f,h\}:=\omega(X_f,X_h)$, where $X_f$ and $X_h$ are the symplectic gradients of $f$ and $g$ respectively. It is well known that $(C^{\infty}(M),\{\cdot,\cdot\})$ is a Lie algebra and that the map
\[
f \mapsto X_f
\]
is an anti-homomorphism of Lie algebras, see \cite{stp}. Moreover, the Poisson bracket is invariant under symplectomorphism, i.e. given a symplectomorphsim $\phi:M\lra M$, of $(M,\omega)$, for all $f,h\in C^\infty(M)$, we have $\{\phi^* f,\phi^* f \}=\phi^* \{f,g\}$; see for instance \cite{Ca}. More in general, if $\phi: M\to N$ is a smooth map between smooth manifolds equipped with Posson brackets $\{\cdot,\cdot\}_M$, $\{\cdot,\cdot\}_N$, respectively, we say that $\phi$ is a \emph{Poisson map} if for all $f,h\in C^\infty(N)$, we have $\phi^* \{f,g\}_N=\{\phi^* f,\phi^* f \}_M$.


Let now $G$ be a Lie group acting by symplectomorphisms on  $(M,\omega)$, i.e. $g^*\omega=\omega$ for all $g\in G$.  From now on, we always assume that $G$ acts properly on $M$ and that $M/G$ is connected.  

\begin{defini}
We say that the $G$-action on $M$ is \textit{coisotropic} if there exists an open dense subset $U$ of $M$ such that for every $x\in U$, $G(x)$ is a coisotropic submanifold.
\end{defini}

On the other hand, by the properties of symplectomorphism that we have previously recalled, we note that the algebra of the $G$-invariant smooth functions on $M$, that we denote as $(C^{\infty} (M)^G,\{\cdot,\cdot \})$, is a Lie subalgebra of $(C^{\infty} (M),\{\cdot,\cdot \})$. We can now give the following definition.

\begin{defini}
A symplectic $G$-action on $(M,\omega)$ is called \emph{multiplicity free} and
$M$ is called a \textit{$G$-multiplicity free space} if
$( C^{\infty}(M)^G, \{\cdot,\cdot\})$ is a commutative Lie algebra.
\end{defini}
In particular, note that if $K\subset G$ is a subgroup and the $K$-action on $M$ is multiplicity free, then the $G$-action on $M$ is multiplicity free as well.  
We claim that a coisotropic action is multiplicity free. Let $\xi \in \lieg$ and let $\xi_M$ be the fundamental vector field induced by $\xi$ by the $G$-action defined in (\ref{infgen}). Since $G$ preserves $\omega$, the vector field $\xi_M$ is a locally Hamiltonian vector field. 
\begin{lemma}\label{invariant}
Let $f\in C^\infty (M)^G$ and let $x\in M$. Then $X_f(x) \in (T_x G(x) )^\omega$.
\end{lemma}
\begin{proof}
Let $\xi \in \lieg$. Since $f$ is $G$-invariant, it follows that
\[
0=(\mathrm d f )_x (\xi_M (x))=\omega_x(X_f(x),\xi_M(x))
\]
and so $X_f(x)  \in (T_x G(x) )^\omega$.
\end{proof}
\begin{prop} \label{multfree iff coiso}
$M$ is a $G$-multiplicity free space if and only the $G$-action on
$M$ is coisotropic.
\end{prop}
\begin{proof}
Assume that $M$ is a $G$-multiplicity free space. Let $x\in \princ$. By the Tubular Neighborhood Theorem, there exists a $G$-invariant neighborhood $U$ of $G(x)$ which is $G$-equivariantly diffeomorphic to the trivial bundle
$G/G_x \times (T_x G(x))^\perp$. Let $k=\mathrm{cohom}(G,M)=\dim (T_x G(x))^\perp$. Then there exists functions $f_1,\ldots,f_k \in C^{\infty}(U)^G$ with
$\mrm d f_1 \wedge \ldots \wedge \mrm d f_k \neq 0$ such that
\[
G(x)=\{ y \in U : f_1(y)=\ldots=f_k(y)=0 \}.
\]
These functions can be provided, for instance, by the coordinate functions on $(T_x(G(x)))^\bot$ pulled-back to the product $(G/G_x)\times (T_x(G(x)))^\bot$ by the projection on the second factor. 
We can extend $f_1,\ldots,f_k$ to smooth $G$-invariant functions on $M$. The argument is standard. For the sake of completeness we sketch the proof.
Let
$\sigma : G/G_x \times (T_x G(x))^\perp \lra \R$  be the function defined by
\[
\sigma((gG_x,v))=\frac{f(1-\parallel v \parallel^2)}{f(1-\parallel v \parallel^2)+f(\parallel v \parallel^2-\frac{1}{4})},
\]
where
\[
f(t)=\left\{\begin{array}{lc} e^{1/t} & t>0 \\ 0 & t\leq 0.\end{array}\right.
\]
Then $\tilde f_i :=\left\{\begin{array}{lc} \sigma(y) f_i (y) & \mathrm{if\ } y\in U \\ 0 & \mathrm{if\ } y \notin U \end{array}\right.$ is a $G$-invariant smooth function such that $\tilde f_i =f_i$ on a $G$-invariant neighborhood $W\subseteq U$ of $G(x)$, for $i=1,\ldots,k$.
 Since the action is multiplicity free, we have $\{\tilde f_i,\tilde f_j \}=0$ and so $X_{f_i}(x) \in T_x G(
x).$ By Lemma \ref{invariant},  $X_{f_i}(x) \in (T_x G(
x))^\omega,$ for $i=1, \ldots, k$. Hence, keeping in mind that  $X_{f_1},\ldots
,X_{f_k}$ are linearly independent on $W$, the vector fields $X_{f_1} (x),\ldots
,X_{f_k} (x)$ form a basis of  $(T_x G(x))^\omega$. This proves that $G(x)$ is
coisotropic for any $x\in \princ$. Therefore, by the Principal Orbit Theorem, the $G$-action is coisotropic.

Assume the $G$-action on $M$ is coisotropic. Let $U$ be an open and dense subset of $M$ such that $G(x)$ is coisotropic for any $x\in U$. Let $f,g\in C^{\infty}(M)^G$. By Lemma \ref{invariant} for any  $x\in U$, we have
$X_f (x), X_h (x) \in (T_x G(x) )^\omega \subseteq T_x G(x)$.  This implies  $\{f,h\}(x)=\omega_x (X_f (x) ,X_h (x))=0$ for any $x\in U$. Therefore, keeping in mind that  $U$ is open and dense, we have $\{f,h\}=0$ on $M$.
\end{proof}
\begin{cor}\label{poi}
The $G$-action on $M$ is coisotropic if and only if $G(x)$ is coisotropic for any $x\in \princ$.
\end{cor}

\begin{cor}\label{connected}
The $G$-action on $M$ is coisotropic if and only if the $G^o$-action on $M$ is coisotropic. If $M'$ denotes a connected component of $M$, then the $G$-action on $M$ is coisotropic if and only if the $G^o$-action on $M'$ is coisotropic.
\end{cor}
\begin{proof}
Assume that $G^o$ acts coisotropically on $M$. Since $G^o \subseteq G$,  keeping in mind that $( C^{\infty}(M)^G, \{\cdot,\cdot\})\subseteq ( C^{\infty}(M)^{G^o}, \{\cdot,\cdot\})$, it follows  that $( C^{\infty}(M)^G, \{\cdot,\cdot\})$ is Abelian and so $G$ acts coisotropically on $M$.

Assume that $G$ acts coisotropically on $M$. Let $x\in M$ be such that $G(x)$ is coisotropic. Since $T_x G(x)=T_x G^o (x)$, it follows that $G^o (x)$ is coisotropic as well. Hence, $G^o$ acts coisotropically on $M$.

Finally,  assume that $G^o$ acts coisotropically on $M'$. Then, there exists an open and dense subset $U$ of $M'$ such that $G^o(x)$ is coisotropic for any $x\in G^o$. Let $M''$ be another connected component of $M$. There exists $h\in G$ such that $hM'=M''$. Therefore
\[
G^o (hx)=hG^o(x),
\]
for any $x\in U$. This implies that $GU$ is an open and dense subset of $M$ such that $G^o (y)$ is coisotropic for any $y\in GU$. The vice-versa is easy to check, concluding the proof.
\end{proof}

Under the assumption that $G$ acts properly and by symplectomorphisms on $(M,\omega)$, there exists a $G$-equivariant almost complex structure $J$ such that:
\begin{itemize}
\item $\omega (J\cdot,J\cdot)=\omega(\cdot,\cdot)$;
\item $\mathtt h (\cdot,\cdot):=\omega(\cdot,J\cdot)$ is a Riemannian metric.
\end{itemize}
For a proof of these facts, see \cite[Proposition 19, p.224]{lb}. The collection $(M,\omega,J,\mathtt h)$ is called an \emph{almost K\"ahler manifold}. 
From now on, we always assume that $(M,\omega,J,\mathtt h)$ is an almost K\"ahler manifold, and that the $G$-action preserves all three structures. 

\begin{rem}
Notice that we can express the symplectic orthogonal of a submanifold $P\subset M$ at some point $x\in M$ in terms of $J$ and $\mathtt h$ as $T_x P^\omega=J(T_x P^\perp)$, where $T_x P^\perp$ is the orthogonal complement of $T_xP$ in $T_xM$ relative to $\mathtt h$. 
It follows that a submanifold $P$ is isotropic, coisotropic, Lagrangian or symplectic if for any $x\in P$ we have
$T_x P\subseteq J(T_x P^\perp)$,
$J(T_x P^\perp)\subseteq T_x P$, $J(T_x P^\perp)= T_x P$, $J(T_xP^\perp) \cap T_x P=\{0\}$ respectively.
\end{rem}

Let $H$ be a principal isotropy for the $G$-action on the almost K\"ahler manifold $(M,\omega,J,\mathtt h)$. Let $N\coloneqq N(H)$ denote the normalizer of $H$ in $G$. We denote by $L$ the quotient $N(H)/H$. Then $L$ is a Lie group acting properly on $\core$. In Lemma \ref{core i union of conn.comps} we prove that $\core$ is the union of the connected components of $M^H$ that meet $M_H$. By a Theorem of Koszul, $M^H$ is a disjoint union of embedded submanifolds and $T_x M^H=(T_x M)^H$. Since $J$ is $G$-equivariant, it follows that $(T_x M)^H$ is a $J$-invariant subspace. This means that $(\core,\omega|_{\core},J|_{\core},\mathtt h|_{\core})$ is an almost K\"ahler manifold. Moreover, we have the $L$-action preserves all the three structures and that $\core/L$ is connected.
\begin{prop}
If the $G$-action on $M$ is coisotropic, then $\dim (L) \geq \mathrm{cohom}(G,M)$.
\end{prop}
\begin{proof}
Let $x\in \princ$. By the above corollary, $G(x)$ is coisotropic. Therefore, recalling that the slice representation at $x$ is trivial and that $J$ commutes with the $G$-action, it follows that
$H$ fixes point-wise the subspace $J((T_x G(x))^\perp)\subseteq T_x G(x)$. By the Tubular Neighborhood Theorem, we have $(T_x G(x))^H=T_x N(x)$, where $H=G_x$, and so $\dim N/H \geq \mathrm{cohom}(G,M)$, concluding the proof.
\end{proof}
We recall that a $G$-orbit $G(x)$ is called \emph{regular} if $\mrm{codim }\left( G(x)\right)=\dim (M/G)$. We denote the set of points belonging to a regular orbit as $M_{\text{reg}}$. Note that, by definition, the dimension of the regular orbits equals the dimension of the principal orbits. By a continuity argument, we can prove that if the principal orbits are coisotropic also the regular ones are.
\begin{prop}
If the $G$-action on $M$ is coisotropic then $G(y)$ is coisotropic for every $y\in M_{\mrm{reg}}$.
\end{prop}
\begin{proof}
Let again $H$ be a principal isotropy. By Principal Orbit Theorem we may assume that $H\subseteq G_y$. Let $\epsilon>0$ be such that $\mathcal N_y=\exp_y (B_\epsilon \cap ( T_y G( y ))^\perp)$ is a slice at $y$. Then there exists $v\in B_\epsilon \cap ( T_y G(y))^\perp$ such that $(G_x)_v=H$. The geodesic
\[
\gamma: [0,1] \lra M, \qquad t \mapsto \exp_y(tv),
\]
satisfies $G_{\gamma(t)}=H$ for $t\in (0,1]$. By the above Corollary, the orbit $G(\gamma(t))$ is coisotropic for any $t\in (0,1]$. Split $\mathfrak g=\mathfrak h\oplus \mathfrak m$. Let $\xi_1,\ldots,\xi_k$ be a basis of $\mathfrak m$. Then $\xi_1(\gamma(t)),\ldots,\xi_k(\gamma(t))$ is a basis of $T_{\gamma(t)} G (\gamma (t))$ for any $t\in [0,1]$. Using a system of coordinates at $y$ and the Gram-Schmidt procedure, there exists $t_o >0$ and smooth vector fields $\nu_1 (t),\ldots,\nu_{n-k}(t)$ along $\gamma_{|_{[0,t_o]}}$ such that $\nu_1 (t),\ldots,\nu_{n-k}(t)$ is a basis of $(T_{\gamma(t)} G(\gamma (t)) )^\perp$ for any $t\in [0,t_o]$. Since $\dim G(y)= \dim G(\gamma (t))$, keeping in mind that $(T_y G( y))^\omega=J((T_y G(y))^\perp)$, it follows that $\dim (T_y G(y))^\omega= \dim (T_{\gamma(t)} G(\gamma (t)))^\omega$ for any $t\in [0,t_o]$.

Let $v\in (T_y G(y))^\omega$. There exists $\alpha_1,\ldots,\alpha_{n-k} \in \R$ such that
\[
v=\alpha_1 J(\nu_1 (0))+\cdots+\alpha_{n-k} J(\nu_{n-k}(0)).
\]
Hence
\[
v=\lim_{t\mapsto 0} \alpha_1 J(\nu_1 (t))+\cdots+\alpha_{n-k} J(\nu_{n-k}(t)).
\]
Therefore
\[
\mathtt h_y (v,\nu_s(0))=\lim_{t \mapsto 0} \mathtt h_{\gamma(t)} (\alpha_1 J(\nu_1 (t))+\cdots+\alpha_{n-k} J(\nu_{n-k}(t))), \nu_s (t))=0,
\]
for $s=1,\ldots,n-k$,
due to the fact that $G(\gamma(t))$ is coisotropic for any $t>0$. Hence $v\in T_y G(y)$ and the result follows.
\end{proof}

We are ready to prove the reduction principle for coisotropic actions.

\begin{lemma}
If $x\in (\core)_{\mrm{princ}}$, then $L(x)$ is coisotropic in $\core$ if and only if $G(x)$ is coisotropic in $M$.
\end{lemma}
\begin{proof}
Let $x\in (\core)_{\mrm{princ}}$. As a consequence of the Tubular Neighborhood Theorem, keeping in mind that the slice representation is trivial, a $G$-invariant neighborhood $U$ of $G(x)$ is $G$-equivariant diffeomorphic to the trivial bundle
\[
G/H \times (T_{x} G(x))^{\perp}.
\]
Therefore,
\[
\core \cap U=L\times (T_{x} G(x))^{\perp},
\]
which implies $T_x L(x)= (T_x G(x) )^H$ and $(T_x L(x) )^{\perp}\cap (T_x(\core))=(T_x G(x) )^{\perp}$. 
In particular, $(T_x G(x) )^\omega=J(( T_x G(x)
)^{\perp})$. Hence, keeping in mind that $J$ is a $G$-equivariant almost complex structure and the slice representation is trivial, we get
$(T_x G(x))^{\omega}\subset (T_x M )^H$. Therefore,
\[
J( (T_x G(x) )^{\perp }) \subset T_x G(x)
\Leftrightarrow J|_{\core}( (T_x L(x) )^{\perp }) \subset (T_x G(
x) )^H =T_x L(x),
\]
and so
$G(x)$ is coisotropic in $M$ if and only if $L(
x)$ is  coisotropic in $\core$.
\end{proof}

\coisotropicreduction*

\begin{proof}
By Theorem \ref{reduction-global}, we have  $G(\core)_{\mrm{princ}}=\princ$. By the above Lemma, if $x\in (\core)_{\mrm{princ}}$, then $L(x)$ is coisotropic in $\core$ if and only if $G(x)$ is coisotropic in $M$.  By Corollary \ref{poi}, it follows that $L$ acts coisotropically on $\core$ if and only if $G$ acts coisotropically on $M$.
\end{proof}

We can also give an alternative proof of the Reduction Principle for Coisotropic actions, using Corollary \ref{ChevCore} and the following lemma.

\begin{lemma}\label{poissrestr}
    Let $G$ act on $(M,\omega)$ by symplectomorphisms, let $H$ be a principal isotropy and let $\core$ be the core relative to $H$. Then the inclusion $\iota:\core\hookrightarrow M$ satisfies $$\iota^*\{f,g\}=\{\iota^*f,\iota^*g\}_{\core}\quad\forall f,g\in C^\infty(M)^G\ ,$$
   where $\{\cdot,\cdot\}_{\core}$ is the restriction of the Poisson bracket to $\core$.
\end{lemma}
The same lemma is stated for $M_H$ instead of $\core$ in [\cite{OR}, Proposition 4.2.9 pages 139-140], however, the proof is the same. We report it anyway, for completeness.
    \begin{proof}
        Take $f\in{C}^\infty(M)^G$, then if $x\in \core$, $X_f(x)\in T_x(\core)$.
     In fact, for all $v\in T_xM$ and $h\in H$, we have  $\omega_x((\text{d}h)_x(X_f(x)),v)=\omega_x(X_f(x),(\text{d}h)^{-1}_x(v))=(\text{d} f)_x((\text{d}h)^{-1}_x(v))=(\text{d} f)_x(v)=\omega_x(X_f(x),v)$.
    Now, $\forall x\in \core$ and $v\in T_x(\core)$ we have 
    $$(\iota^*\text{d} f)_x(v)=(\text{d}(\iota^*f))_x(v)=(\text{d} f)_x((\text{d} \iota)_x(v))=\omega_x(X_f(x),(\text{d}\iota)_x(v)) .$$
     On the other hand, we can write also
    $$(\iota^*\text{d} f)_x(v)\!=\!(\text{d}(\iota^*f))_x(v)\!=\!(\iota^*\omega)_x(X_{\iota^*f}(x),v)\!=\!\omega_x((\text{d}\iota)_x(X_{\iota^*f}(x)),(\text{d}\iota)_x(v))\ .$$
    Therefore, $\omega_x(X_f(x)-(\text{d}\iota)_x(X_{\iota^*f}(x)),w)=0$, $\forall w\in T_x(\core)$. Since $\core$ is symplectic, that is, $\forall x\in \core$, $T_xM=T_x(\core)\oplus (T_x(\core))^\omega$, $\omega_x(X_f(x)-(\text{d}\iota)_x(X_{\iota^*f}(x)),v)=0$ $\forall v\in T_xM$, and therefore, $$X_f(x)=(\text{d}\iota)_x(X_{\iota^*f}(x))\ .$$
 Take now $f,g\in {C}^\infty(M)^G$ and therefore $\iota^*f,\iota^*g\in {C}^\infty(\core)^{N(H)}$. 
We have $\forall x\in \core$,
$$\{\iota^*f,\iota^*g\}_{\core}(x)=(\iota^*\omega)_x(X_{\iota^*f}(x),X_{\iota^*g}(x))=\quad$$ $$=\omega_x((\text{d}\iota)_x(X_{\iota^*f}(x)),(\text{d}\iota)_x(X_{\iota^*g}(x)))=\qquad\ \ \ \ $$ $$=\omega_x(X_f(x),X_g(x))=\{f,g\}(x)=\iota^*\{f,g\}(x)\ .$$
concluding the proof.
    \end{proof}
The proof of the Reduction Principle for Coisotropic actions then easily follows. In fact, the action of $G$ is coisotropic if and only if $(C^\infty(M)^G,\{\cdot,\cdot\})$ is Abelian. By Corollary \ref{ChevCore} and Lemma \ref{poissrestr}, this is equivalent to $(C^\infty(\core)^{L},\{\cdot,\cdot\})$ being Abelian, which in turn is equivalent to the $L$-action on $\core$ being coisotropic. 

Furthermore, we point out that the reduction principle also holds for a connected component of $\core$. In fact, we have:

\begin{thm}
Let $F$ be a connected component of $\core$. Then $L^o$ acts coisotropically on $F$ if and only if $G$ acts coisotropically on $M$.
\end{thm}
\begin{proof}
By Corollary \ref{connected}  we get that $L^o$ acts coisotropically on $F$ if and only if $L$ acts coisotropically on $\core$. Applying the Reduction Theorem for Coisotropic actions, we get that $L^o$ acts coisotropically on $F$ if and only if $G$ acts coisotropically on $M$.
\end{proof}
Summing up, we have proved the following result.

\equivalencethmcoisotropic*

\section{Infinitesimally almost homogeneous actions}\label{sect inf almost hom}
Let $M$ be a compact connected complex manifold and let $H$ be a closed
complex subgroup of the group of holomorphic automorphisms $\mathrm{Aut}(M)$. If $H$  has an
open orbit $\Omega$ in $M$ then $M$ is called $H$-\emph{almost homogeneous}. The open orbit $\Omega$ is a dense, open and connected complex submanifold
of $M$ and its complement $E: =M\setminus\{\Omega\}$ is a proper analytic subset of $M$, possibly
empty \cite{hs}. Such manifolds arise quite
naturally in many different settings, see for instance \cite{hs,HW}. 
In our setting, $(M,\omega,J,\mathtt h)$ is an almost K\"ahler manifold and a Lie group $G$ acts properly on $M$ preserving all the three structures and $M/G$ is connected. Since $G^\C$ does not act in general on $M$, we give the following definition.
\begin{defini}
The $G$-action on $M$ is called \emph{infinitesimally almost homogeneous} if there exists $x\in M$ such that $T_x M =T_x G(x)+J(T_x G(x))$. 
\end{defini}
\noindent Note that since $T_x M =T_x G(x) \oplus T_x G(x)^\perp$, it follows that condition $T_x M =T_x G(x)+J(T_x G(x))$ holds if and only if $T_x G(x)^\perp \subseteq T_x G(x)+J(T_x G(x))$.
\begin{rem}
If the $G$-action on $M$ is coisotropic, then it is infitesimally almost homogeneous. In fact, let $x\in \princ$. Then $G(x)$ is coisotropic, and so $J(T_x G(x)^\perp ) \subseteq T_x G(x)$. This implies $T_x G(x)^\perp \subseteq J(T_x G(x))$ 
concluding the proof.
\end{rem}

The following lemma is easy to check.
\begin{lemma}\label{connected-almost}
The $G$-action on $M$ is infinitesimally almost homogeneous if and only if the $G^o$-action on $M$ is infinitesimally almost homogeneous. If $M'$ is a connected component of $M$ then the $G$-action on $M$ is infinitesimally almost homogeneous if and only if the $G^o$-action on $M'$ is infinitesimally almost homogeneous.
\end{lemma}

Since the condition $T_x G(x)+J(T_x G(x))=T_xM$ is an open condition, we get the following result.
\begin{prop}\label{infhom}
The $G$-action on $M$ is almost homogeneous if and only if there exists $x\in \princ$ such that $T_x M =T_x G(x)+J(T_x G(x))$. 
\end{prop}
We are ready to prove a reduction principle for almost homogeneous actions.
\infhomreduction*
\begin{proof}
 It suffices to prove the theorem for the $N$-action on $\core$. Let $x\in \left(\core\right)_\mrm{princ}$. Then $T_x N(x)^\perp \cap T_x \, \core =T_x G(x)^\perp$. Since $(T_x G(x))^H=T_x N(x)$ we have
\[
(T_x G(x) + J(T_x G(x))^H=T_x N(x) + J(T_x N(x)).
\]
The subspace $T_x G(x)^\perp$ is fixed point-wise by $H$. This implies
$(T_x G(x))^\perp \subseteq T_x G(x) + J(T_x G(x))$ if and only if $(T_x G(x))^\perp \subseteq T_x N(x) + J(T_x N(x))$ and so if and only if the orthogonal complement of $T_xN(x)$ in $T_x(\core)$ is contained in $T_xN(x)+J\left(T_xN(x)\right)$.
\end{proof}
By Lemma \ref{connected-almost} also the following result holds.
\begin{thm}
If $F$ is a connected component of the $\core$, then the $L^o$-action on $F$ is infinitesimally almost homogeneous if and only if the $G$-action on $M$ is infinitesimally almost homogeneous.
\end{thm}

\section{Coisotropic Hamiltonian Actions}\label{sect hamiltonian}
Let  $G$ be a Lie group acting properly and symplectically on $(M,\omega)$ with $M/G$ connected.  Let $\xi\in \mathfrak g$ and $\xi_M$ be its associated fundamental vector field. Recall that $\xi_M$ is locally Hamiltonian, that is, $\iota_{\xi_M}\omega$ is closed.
\begin{defini}
    We say that the $G$-action on $M$ is \emph{Hamiltonian} if there exists a map $\mu:M\to \mathfrak g^*$, where $\mathfrak g^*$ is the dual of the Lie algebra $\mathfrak g$ of $G$, called the \textit{momentum map}, satisfying:
    \begin{itemize}
        \item For any $\xi\in \mathfrak g$, let $\mu^\xi:M\to \R$, $p\mapsto \langle \mu (p),\xi\rangle$,  where $\langle \cdot ,\cdot \rangle$ is the dual pairing; then $\iota_{\xi_M}\omega=d\mu^\xi$;
        
        \item $\mu$ is $G$-equivariant with respect to the coadjoint action of $G$ on $\mathfrak g^*$; i.e.
        $    \mu (gx)= g \mu(x)$ for any $g\in G$ and any $x\in M$, where $g \mu(x)=\bigr(\text{Ad}(g^{-1})\bigr)^*\mu(x)$.
        
    \end{itemize}

\end{defini}

\begin{lemma}\label{lemma prop's of mom map}
    Let $x\in M$. Then the following facts hold true:
    \begin{enumerate}
        \item $\ker d\mu_x=(T_xG(x))^\omega$;
        \item $\text{im } (d\mu)_x=\mathfrak g_x^\perp=\{\phi \in \mathfrak g^*\ :\  \phi(v)=0, \text{ for all } v\in \mathfrak g_x\}$. Thus $(d\mu)_x$ is surjective if and only if $\mathfrak g_x=\{0\}$;
        \item $\omega\vert_{G(x)}=\bigr(\mu^*\Omega\bigr)\bigr \vert_{G(x)}$, where $\Omega$ is the Kirillov-Kostant-Souriau symplectic form defined on the coadjoint $G$-orbit through $\mu(x)$;
        \item $(T_xG(x))\cap\bigr( T_xG(x)\bigr)^\omega=T_xG_{\mu(x)}(x) $;
        \item if $x\in \princ$, then 
        \begin{align*}
            \dim (G_{\mu(x)})-\dim (G_x)\leq \text{cohom}(G,M).
        \end{align*}

        \end{enumerate}
\end{lemma}
\begin{proof}
$v\in \mathrm{ker}\, ({\rm d} \mu)_x $ if and only if for any $\xi \in \lieg$ we have $\mrm (d \mu^\xi)_x=0$. Since $(\mrm d \mu^\xi )_x=(i_{\xi_M}\omega)(x)$, it follows that $v\in \mathrm{ker}\, ({\rm d} \mu)_x $ if and only if $v\in (T_x G(x))^\omega$.

Let $\xi \in \lieg_x$ and let $v\in T_x M$. Then $(\mathrm{d} \mu^\xi)_x (v)=(i_{\xi_M (x)}\omega) (v)=0$. Hence  $\mrm{im}\,(\mrm d \mu)_x\subset  \mf g^o_x$. On the other hand
\[
\dim (\mathrm{im}\, (\mrm d \mu)_x) =\dim (M) -\dim (\mathrm{ker}\, (\mathrm{d} \mu)_x)=\dim (M) - (\dim (M) - \dim (G) + \dim (G_x))=\dim ( \lieg_x^o),
\]
and so item \emph{(2)} holds by dimensional reason.

Let $\xi,\nu \in \lieg$. Then
\[
\mu(\exp(t\xi)x) (\nu)=\mathrm{Ad}(\exp(-t\xi))^*(\mu(x))(\nu),
\]
and so
\[
(\mathrm d \mu^\nu)_x (\xi_M (x)) =\omega(\nu_M (x),\xi_M (x))=\mu(x)([\nu,\xi])=:\Omega(\nu_{G( \mu(x))} (x),\xi_{G(\mu(x))} (x)).
\]
Therefore, the restriction of the symplectic form $\omega$ to $G(x)$ coincides with the pull-back by the momentum map $\mu$ of the Kirillov-Kostant-Souriau symplectic form  defined on the coadjoint orbit through $\mu(x)$ \cite{stp}. In particular,
\[
(T_x G(x) ) \cap (T_x G(x)  )^{\omega}= T_x G_{\mu(x)} (x),
\]
and hence, $\dim (G_{\mu(x)})- \dim (G_x) = \dim \bigl( T_x G(x) \,\cap\, (T_x G(x) )^\omega \bigr)$. Item \emph{(5)} easily follows.
\end{proof}
The momentum map gives a criterion for a $G$-orbit being symplectic or coisotropic.
\begin{prop}\label{prop characterization mom map symp coiso}

    Let $x\in M$. The following assertions hold true:
    \begin{enumerate}
        \item $G(x)$ is symplectic if and only if $(G_x)^0=(G_{\mu (x)})^0$ and so if and only if $\mu: G(x)\to G(\mu(x))$ is a covering map;
        \item If $x\in \princ$, then $G(x)$ is coisotropic if and only if $\dim (G_{\mu(x)})-\dim (G_x)=\text{cohom}(G,M)$.
    \end{enumerate}
    
\end{prop}
\begin{proof}
    By the Lemma \ref{lemma prop's of mom map} \emph{(4)}, we have that $G(x)$ is symplectic if and only if $T_xG_{\mu(x)}(x)=0$ and so if and only if $\mathfrak g_x=\mathfrak g_{\mu(x)}$. This proves item \emph{(1)}

    Since $\dim ((T_xG(x))^\omega)=\text{cohom}(G,M)$, we have that item \emph{(2)} follows from Lemma \ref{lemma prop's of mom map} \emph{(5)}.
\end{proof}
We now recall a theorem by Bates and Lerman that will be useful for the reduction principle in the coisotropic Hamiltonian case. 
As usual, let $H$ be a fixed principal isotropy; let $N=N(H)$ its normalizer in $G$, and let $L=N/H$.

\begin{prop}[{{\cite[Theorem 6, page 220]{lb}}}]\label{G hamil implies N(H) hamil}
    If the $G$-action on $M$ is Hamiltonian, then the $L$-action on $\core$ is Hamiltonian.
\end{prop}

From now on, we will assume further that the coadjoint $G$-orbits in $\mathfrak g^*$ are locally closed and hence embedded submanifolds of $\mathfrak g^*$. This condition ensures that the local normal form of the momentum map, which we now describe following \cite{lb}, holds; for further details, we refer to \cite{lb, sl}. 

Let $x\in M$ and $\beta=\mu(x)$. Let $V$ be an $\mathtt{h}$-orthogonal complement of $T_xG_{\beta}(x)$ inside $(T_xG(x))^\omega$. It follows that $V$ is a symplectic subspace of $T_xM$. Keeping in mind that $G_x$ is compact, we may fix an $\text{Ad}(G_x)$-invariant inner product on $\mathfrak g^*$ and note that $\mathfrak g_{\beta}$ and $\mathfrak g_x$ are both $G_x$-invariant. Let $\mathfrak m$ be the orthogonal complement of $\mathfrak g_x$ inside $\mathfrak g_{\beta}$ and let $\mathfrak s$ be the orthogonal complement of $\mathfrak g_{\beta}$ inside $\mathfrak g$.  Thus, we have 
\begin{align} \label{decomp lie alg. coadj}
    \mathfrak g = \mathfrak g_\beta \oplus \mathfrak s=\mathfrak g_x \oplus \mathfrak m \oplus \mathfrak s.
\end{align}

Consider $\mathfrak s x:= \{\xi_M(x)\in T_xM\ :\  \xi \in \mathfrak s\}$. This is also a symplectic subspace of $T_xM$ and we have that $V\cap (\mathfrak s x)=\{0\}$. We therefore have the following. 
\begin{align*}
    T_xM=\bigr(V\oplus \mathfrak s x\bigr)\oplus \bigr(V\oplus \mathfrak s x\bigr)^\omega. 
\end{align*}
Note that $T_xG_\beta(x)$ is a Lagrangian subspace of the symplectic space $\bigr(V\oplus \mathfrak s x\bigr)^\omega$. Since $T_xG_\beta (x)\cong \mathfrak g_\beta /\mathfrak g_x$, we thus get
\begin{align*}
    (T_xM,\omega(x))\cong (V,\omega_x\vert_V)\oplus (\mathfrak s x, (\mu^*\omega)|_{G(x)}) \oplus ((\mathfrak g_\beta /\mathfrak g_x) \oplus (\mathfrak g_\beta /\mathfrak g_x)^*, \omega_{\text{std}})\ ,
\end{align*}
where $\omega_{\text{std}}$ is the canonical symplectic form on the direct sum of a vector space and its dual\footnote{Explicitly, if $W$ is a vector space and $W^*$ is its dual, given $v,w\in W$ and $\alpha,\beta\in W^*$, one has that $\omega_{\text{std}}((v,\alpha),(w,\beta)):=\beta(v)-\alpha(w)$ is a symplectic form on $W\oplus W^*$.}. 
Since $V\oplus (\mathfrak g_\beta /\mathfrak g_x)^*$ gives a $G_x$-invariant complementary space of $T_xG(x)$, we get that there is a $G$-invariant neighborhood of $G(x)$ that is $G$-equivariantly symplectomorphic (see \cite{lb}) to a $G$-invariant neighborhood of the zero section of $Y=G\times_{G_x}\bigr(V\oplus (\mathfrak g_\beta /\mathfrak g_x)^*\bigr) \to G/G_x\cong G(x)$. With this identification, the momentum map $\mu_Y:Y\to \mathfrak g^*$ is given by 

\begin{align*}
    \mu_Y([g,v,m])=g \bigr( \beta + i(\mu_V(v)) + j(m)\bigr),
\end{align*}

where $\mu_V:V\to \mathfrak g_x^*$ is the momentum map of the linear $G_x$-action on the symplectic vector space $(V,\omega_x\vert_V)$. Finally, $i:\mathfrak g_x^*\to \mathfrak g^*$ and $j:\mathfrak g_\beta^*\to \mathfrak g^*$ are the inclusions induced by the decomposition (\ref{decomp lie alg. coadj}). 

Note that if $G(x)$ is symplectic, then $T_xM=(T_xG(x))\oplus \bigr(T_xG(x)\bigr)^\omega$ and thus $\mathfrak g_\beta=\mathfrak g_x$. We have then that $V=\bigr(T_xG(x)\bigr)^\omega$ and thus there is a $G$-invariant neighborhood of $G(x)$ symplectomorphic to a neighborhood of the zero section of $Y=G\times_{G_x}V\to G/G_x$. The tangent $G_x$-action on $(V,\omega_x\vert_V)$ is called the \textit{symplectic slice representation}. Since $\bigr(T_xG(x)\bigr)^\omega=J\bigr(T_xG(x)\bigr)^\perp$ and $J$ is $G$-equivariant, then the symplectic slice representation is equivalent to the normal slice representation. Following \cite{Kac} and Theorem \ref{equivalence}, we say that the symplectic slice representation is multiplicity free if every principal $G_x$-orbit in $V$ is coisotropic.

\begin{prop}\label{G action on Y coiso iff symp slice mult free}
    Assume that $G(x)$ is symplectic, then the $G$-action on $Y$ is coisotropic if and only if the symplectic slice representation is multiplicity free. 
\end{prop}
\begin{proof}
    Assume that the $G$-action on $Y$ is coisotropic. We will show that $G_x(v)$ is coisotropic for any $v\in V_{\mrm{princ}}$. Since $Y_\text{princ}=G V_{\text{princ}}$, this will suffice. Let $v\in V$ such that $y=[e,v]\in Y=G\times_{G_x}V$ is a principal point and thus $G(y)$ is coisotropic. Let $X\in (T_vG_x( v))^\omega=\ker (d\mu_V)_v$. By the local normal form of the momentum map, we have that $$\mu_Y([g,v])=g(\beta + i\mu_V(v)).$$
    Therefore $(d\mu_Y)_{y}([0,X])=(d\mu_V)_v(X)=0$. Hence, $[0,X]\in \ker (d\mu_Y)_{y}=(T_yG( y))^\omega\subset T_yG(y)$. Since $V$ defines a slice for the $G$-action at $x=[e,0]$, we have that $(T_yG(y))\cap (T_yV)=T_yG_x(y)$. It follows that $X\in T_vG_x(v)$. We have thus shown that $(T_vG_x(v))^\omega \subset T_vG_x(v)$.

    Conversely, suppose the slice representation is multiplicity free. Let $v\in V_\text{princ}$. We have that  $y=[e,v]\in Y_\text{princ}$. By hypothesis $G_x(v)$ is coisotropic. By Proposition \ref{prop characterization mom map symp coiso} \emph{(2)}
    \begin{align*}
        \dim \left((G_x)_{\mu_V(v)}\right)-\dim (G_y)=\text{cohom}(G_x,V)=\text{cohom}(G,M).
    \end{align*}
    Since $G_x\subset G_\beta$ by the $G$-equivariance of the momentum map, we get that $(G_x)_{\mu_V(v)}\subset G_\beta$. It follows that 
    \begin{align*}
        \dim (G_\beta) -\dim (G_y) \geq \text{cohom}(G,M).
    \end{align*}
    Since the reversed inequality holds by Lemma \ref{lemma prop's of mom map} \emph{(5)}, we get equality. By Proposition \ref{prop characterization mom map symp coiso} \emph{(2)}, the orbit $G(y)$ is coisotropic. Keeping in mind that $Y_\text{princ}=G V_{\text{princ}}$, the proposition follows.
 \end{proof}

 P\"uttmann \cite{Pu} introduced the homogeinity rank of the action of a compact Lie group $G$ on $M$ as the integer
\begin{align*}
    \text{homrank}(G,M):=\text{rank} (G)-\text{rank} (H)-\text{cohom}(G,M),
\end{align*}
where $H$ is a principal isotropy group of the $G$-action on $M$. Huckleberry and Wurzbacher proved that a Hamiltonian action of a compact Lie group has vanishing homogeneity rank  if and only if the principal orbits are coisotropic \cite[Theorem 3, page 267]{HW}. We also point out that the classification of the real irreducible representations of compact Lie groups with vanishing homogeneity rank is given in \cite{pg}.

\begin{cor}
    Let $x\in M$ be such that $G(x)$ is symplectic and let $G_x\to O((T_xG(x))^\perp) $ be the normal slice representation. Then the $G$-action on $Y=G\times_{G_x} V$ is coisotropic if and only if $\text{homrank}(G_x,(T_xG(x))^\perp)=0$. If $G$ is compact, then $G$ acts coisotropically on $M$ if and only if $\text{homrank}(G_x,(T_xG(x))^\perp)=0$.
\end{cor}
\begin{proof}
    By Proposition \ref{G action on Y coiso iff symp slice mult free}, the $G$-action on $Y$ is coisotropic if and only if the symplectic slice representation is multiplicity free. By the theorem \cite[Theorem 3, page 267]{HW} stated above, this is equivalent to the vanishing of the homogeneity rank of the $G_x$-action on $J(T_xG(x))^\perp$. Since $J$ is $G$-equivariant, the first claim follows.

    Now assume that $G$ is compact. Since the stabilizers of the coadjoint action are connected and using Proposition \ref{prop characterization mom map symp coiso}, keeping in mind that $G(x)$ is symplectic, we get that $G_x=G_\beta$, where $\beta=\mu(x)$. Hence, $G_x$ contains a maximal torus of $G$. A straightforward computation shows that 
    \begin{align*}
        \text{homrank}(G,M)=\text{homrank}(G_x,(T_xG(x))^\perp).
    \end{align*}
    The theorem of Huckleberry and Wurzbacher \cite[Theorem 3, page 267]{HW}, allows us to complete the proof. 
\end{proof}

We will now discuss the Hamiltonian reduction of the $G$-action on $(M,\omega,J,\mathtt h)$. We will follow \cite{acg, lb,sl}. Let $\mu:M\to \mathfrak g^*$ be the momentum map and $\alpha\in \mathfrak g^*$. We define the corresponding reduced space
\begin{align*}
    M_\alpha = \mu^{-1}(G(\alpha)) /G,
\end{align*}
to be the topological quotient of the subset $\mu^{-1}(G(\alpha))$ of $M$ by the action of $G$. The reduced space $M_\alpha$ is a stratified space where the strata, called \textit{pieces}, are symplectic manifolds. The space $M_\alpha$ can further be endowed with a Poisson struture $\bigr(C^\infty(M_\alpha),\{\cdot,\cdot\}\bigr)$, where 
\begin{align}
    C^\infty(M_\alpha)\coloneqq C^\infty(M)^G\Bigr\vert_{\mu^{-1}(G(\alpha))};
\end{align}
that is, the smooth functions on $M_\alpha$ are the restrictions of $G$-invariant smooth functions on $M$ to the subset $\mu^{-1}(G(\alpha))$. The Poisson bracket $\{\cdot,\cdot \}$ on $M_\alpha$ is induced by the Poisson bracket on $C^\infty(M)$. The pieces of $M_\alpha$ are given as follows. Let  $H\subset G$ be an isotropy group of the $G$-action on $M$. Then, we recall that any connected component $C$  of $ M_{(H)}\coloneqq \{x\in M\ :\  G_x \text{ is conjugate to } H\}
$ is a smooth submanifold of $M$. 
Since we assume that the coajoint orbits are locally closed, it follows by the Symplectic Slice Theorem \cite{OR}, that the subset $\mu^{-1}(G(\alpha))\cap C$ are submanifolds of $M$ such that the topological quotient
\begin{align*}
    (M_\alpha)_{C}\coloneqq \Bigr(\mu^{-1}(G(\alpha)) \cap C \Bigr) \bigr/ G
\end{align*}
has the structure of a smooth symplectic manifold (see \cite{lb}) satisfying:
\begin{itemize}
    \item The restriction of the projection $\pi:M\to M/G$ to $\mu^{-1}(G(\alpha) )\cap C\to (M_\alpha)_C$ is a smooth map;
    \item The inclusion $(M_\alpha)_C\hookrightarrow M_\alpha$ is Poisson map. 
\end{itemize}

\begin{lemma} \label{lemma coisotrop iff red spac discrete}
    Let $(M,\omega)$ be a symplectic manifold and let $G$ be a Lie group acting in a Hamiltonian fashion on $M$.  Then the $G$-action on $M$ is coisotropic if and only if $M_\alpha$ is discrete for every $\alpha\in\mathfrak g ^*$. 
\end{lemma}
\begin{proof}
    Assume that the $G$-action on $M$ is coisotropic. By Theorem \ref{equivalence}, the Poisson structure $\bigr(C^\infty(M)^G,\{\cdot,\cdot\}\bigr)$ is abelian. Hence, also the Poisson structure of $M_\alpha$ is abelian for any $\alpha\in \mathfrak g^*$. It clearly suffices to show that each piece of $M_\alpha$ is discrete. Let $S$ be a piece of $M_\alpha$. We have that the inclusion $S\hookrightarrow M_\alpha$ is a Poisson map, and thus the Poisson structure of $S$ is also abelian. Since the Poisson structure on $S$ comes from its symplectitc structure, it follows that $S$ is zero dimensional and consists of a discrete space of points.

    Conversely, suppose that all reduced spaces are discrete. By the Theorem \ref{equivalence}, it suffices to show that $\bigr(C^\infty(M)^G,\{\cdot,\cdot\}\bigr)$ is abelian. Let $f,g\in C^\infty(M)^G$. Let $x\in M$ and $\alpha=\mu(x)$. Let $C$ be a connected component of $M_{(G_x)}$, so that $C\cap \mu^{-1}(G(\alpha))$ contains $x$ and $(M_\alpha)_C$ is a piece of $M_\alpha$. Since $M_\alpha$ is discrete, so is $(M_\alpha)_C$. Therefore $\{f,g\}(x)=0$. Since $x$ was chosen arbitrarily, $\{f,g\}\equiv 0$.
\end{proof}

\begin{cor}
    Let $(M,\omega)$ be a symplectic manifold and let $G$ be a Lie group acting in a Hamiltonian fashion on $M$.  Then the $G$-action on $M$ is coisotropic if and only if $(\core)_\alpha$ is discrete for every $\alpha\in\mathfrak l ^*$, where $\mathfrak l$ is the Lie algebra of $L$
\end{cor}

\begin{proof}
    Since the $G$-action on $M$ is Hamiltonian, then by Proposition \ref{G hamil implies N(H) hamil} so is the $L$-action on $\core$. Theorem \ref{equivalence} says that the $G$-action on $M$ is coisotropic if and only if the $L$-action on $\core$ is coisotropic. Since the coadjoint orbits of $L$ on $\mathfrak l^*$ are locally closed, then by Lemma \ref{lemma coisotrop iff red spac discrete} this is equivalent to $(\core)_\alpha$ being discrete for all $\alpha\in \mathfrak l^*$.
\end{proof}

Summing up, we have proved the following result.

\hamiltoniancoisoreduction*
\begin{exe}
    Consider the diagonal action of $U(n)$ on $M=\underbrace{\C^ n\oplus\cdots \oplus \C^ n}_{m\text{ times }}$ with $n\geq 2$. Consider the functions $f,g:M\to \R$ given by
    \begin{align*}
        f(z_1,\dots,z_m)&=\Re(z_1^ *z_2), & g(z_1,\dots,z_m)&=\Re(z_2^*z_3);
    \end{align*}
    where $\Re(\cdot)$ denotes the real part. The functions $f$ and $g$ are both $U(n)$-invariant. Moreover, it is easy to see that $\{f,g\}\neq 0$. Hence, for $m\geq 3$, we have that this action is not coisotropic. If $m=1$, we have that the principal $U(n)$-orbits in $\C^ n$ are spheres. Since every codimension $1$ submanifold is coisotropic, we have that the $U(n)$ action on $\C^ n$ is coisotropic. 
    
    Now consider the action of $U(n)$ on $\C^ n\oplus \C^ n$. The principal isotropy is $U(n-2)$. Then $\bigr(\C^ n\oplus \C^ n\bigr)^ {U(n-2)}=\C^2\oplus \C^ 2$ and the $(N(U(n-2))/U(n-2))$-action on $\C^ 2\oplus \C^ 2$ is equivalent to the diagonal $U(2)$ action on $\C^2\oplus \C^ 2$. The cohomogeinity of the latter is $4$ and the principal isotropy is trivial. The momentum map is given by
    \begin{align*}
        \mu(z,w)=-\frac{\mathbf i}{2} (zz^*+ww^*).
    \end{align*}
    Let $\{e_1,e_2\}$ be the standard basis of $\C^2$. Then the orbit of $x_0=(e_1,e_1+e_2)\in \C^2\oplus \C^2$ is principal. Moreover,
    \begin{align*}
        \mu(x_0)=\frac{\mathbf i}{2} \begin{pmatrix}
            2 & 1\\
            1 & 0
        \end{pmatrix}.
    \end{align*}
    Hence $\dim U(2)_{\mu(x_0)}<4$. By Theorem \ref{equiham}, the $U(n)$ action on $\C^ n\oplus \C^ n$ is not coisotropic. Summing up, we have proved that the diagonal action of $U(n)$ on $\underbrace{\C^{n}\oplus\cdots \oplus \C^ n}_{m}$ is coisotropic if and only if $m=1$. This example shows that if a principal $G$-orbit is coisotropic then the $G$-action is not in general coisotropic. Indeed, let  $(e_1,e_2)\in \C^2\oplus \C^2$. Then
    \begin{align*}
        \mu(e_1,e_2)=\frac{\mathbf i}{2}\begin{pmatrix}
            1 & 0 \\
            0 & 1 
        \end{pmatrix},
    \end{align*}
    and so $\dim (\mathrm{U}(2)_{\mu(e_1,e_2)})=4=\text{cohom}(U(2),\C^4)$. By Proposition \ref{prop characterization mom map symp coiso} the ${U}(2)$-orbit throught $(e_1,e_2)$ is coisotropic. Hence the
    $U(n)((e_1,e_2))$ is coisotropic but the  $U(n)$-action on $\C^n\oplus \C^n$ is not coisotropic. Note that the diagonal action of $U(n)$ on $\C^n\oplus \C^n$ is infinitesimally almost homogeneous by Theorem \ref{red infhom} and Corollary \ref{infhom}.
\end{exe}

\begin{lemma} (Restriction Lemma).
    Let $C$ be a closed $G$-invariant symplectic submanifold of $M$. If the $G$-action on $M$ is coisotropic, then so is the $G$-action on $C$.
\end{lemma}
\begin{proof}
    The restriction of the momentum map $\mu$ to $C\to \mathfrak g^*$ is a $G$-equivariant momentum map for the $G$-action on $C$. For $\alpha\in \mathfrak g^*$, we have that $C_\alpha\subset M_\alpha$ and so $C_\alpha$ is discrete. By Theorem \ref{equiham}, it follows that the $G$-action on $C$ is coisotropic.
\end{proof}

\begin{cor}
    Let $G$ be an abelian Lie group acting properly, coisotropically and in a Hamiltonian fashion on $(M,\omega)$. Let $H\subset G$ be a compact subgroup. Then $M^H$ is $G$-invariant and the $G$-action on $M^H$ is coisotropic.
\end{cor}
\begin{proof}
    $M^H$ is $G$-invariant because $G$ is abelian; it is also symplectic because $J$ is $G$-invariant and $J(T_xM^H)\subset T_xM^H$.
\end{proof}
\begin{cor}
    If $G$ is a compact Lie group acting coisotropically on $(M,\omega)$, then $M^G$ must be discrete.
\end{cor}

We refer to Bates and Lerman \cite[Theorem 16, page 220]{lb} for the following. Let $H$ be a compact subgroup of $G$ and let $L=N(H)/H$. A connected component $F$ of $M_H=\{x\in M\ :\  G_x=H\}$ is a symplectic submanifold of $M$ and the $L$-action on $F$ is Hamiltonian, proper and free. Let $x\in F$ and $\alpha=\mu(x)$, Bates and Lerman further showed that there is some $\alpha_0\in \mathfrak l^*$ such that the reduced space $F_{\alpha_0}$ is isomorphic to the symplectic piece $$(M_{\alpha})_F=\bigr(\mu^{-1}(G(\alpha)) \cap F \bigr)/G.$$ 
Using Lemma \ref{lemma coisotrop iff red spac discrete} we obtain the following result.
\begin{prop}
    Let $H$ be a compact subgroup of $G$. If $G$ acts coisotropically on $M$, then $L$ acts coisotropically on $F$.
\end{prop}

Assume now that $G$ is compact and acts coisotropically on $M$. By the theorem of Huckleberry and Wurzbacher \cite{HW}, $\text{homorank}(G,M)=0$. Let $H$ be a principal isotropy and $N=N(G)$ its normalizer in $G$. By Theorem \ref{equiham}, the $N$-action on a connected component $F$ of $\core$ is coisotropic and hence $\text{homorank}(N,F)=0$. Since $\text{cohom}(G,M)=\text{cohom}(N^0,F)$, we thus obtain that $\text{rank}(G)=\text{rank}(N^0)$. This proves the following.
\begin{prop}
    Let $G$ be a compact Lie group acting coisotropically on $(M,\omega)$. Let $H$ be a principal isotropy. Then $N(H)^0$ contains a maximal torus of $G$.
\end{prop}

Let us now consider a Lie group $G=G_1\times G_2$, where $G_1$ and $G_2$ are closed Lie subgroups. We assume that the coadjoint orbits of $G_1$ and of $G_2$ are locally closed in $\mathfrak g_1^*$ and $\mathfrak g_2^*$, respectively. Obviously $\mathfrak g^*=\mathfrak g_1^*\oplus \mathfrak g_2^*$, and the momentum map $\mu:M\to \mathfrak g^*$ decomposes as $\mu=\mu_1+\mu_2$, where $\mu_i$ is the corresponding momentum map for the $G_i$-action on $M$ for $i=1,2$.

Let $\alpha_2\in\mathfrak g_2 ^*$. The $G_1$ action on the pieces of $M_{\alpha_2}$ is symplectic and Hamiltonian. Thus we may obtain a continuous map 
\begin{align*}
\mu_{12}:M_{\alpha_2}\to \mathfrak g_1^*, \quad [x]\mapsto \mu_1(x).  
\end{align*}
which we also call momentum map, whose restriction to each piece gives the momentum map correspoding to the respective $G_1$-action.

\begin{defini}
    We say that the $G_1$-action on $M_{\alpha_2}$ is \textit{multiplicity free} if the ring of $G_1$-invariant smooth functions on $M_{\alpha_2}$ forms an abelian Poisson subalgebra of $C^\infty(M_{\alpha_2})$
\end{defini}

For any $\alpha_1\in \mathfrak g^*_1$, we consider the reduced space with respect to the $G_1$-action on $M_{\alpha_2}$ defined as
\begin{align*}
    (M_{\alpha_2})_{\alpha_1}=\mu_{12}^{-1}(G_1(\alpha_1))/G
\end{align*}
 We make similar remarks and definitions for the $G_2$-action on $M_{\alpha_1}$. 
We point out that the $G_1$-action on $M_{\alpha_2}$ preserves the pieces of $M_{\alpha_2}$. We finally note that if we write $\alpha=\alpha_1+\alpha_2\in \mathfrak g^*=\mathfrak g_1^*\oplus \mathfrak g_2^*$, then $ (M_{\alpha_2})_{\alpha_1}$ are both $ (M_{\alpha_1})_{\alpha_2}$ homeomorphic to $M_\alpha$ via the natural map $[[x]]\mapsto [x]$.

\begin{thm}
    Let $(M,\omega)$ be a symplectic manifold and let $G=G_1\times G_2$ be a Lie group acting in a Hamiltonian fashion on $M$. Assume that the coadjoint orbits of $G_1$ and $G_2$ are respectively locally closed submanifolds. The following are equivalent
    \begin{enumerate}
        \item The $G$-action on $M$ is multiplicity free;
        \item For every $\alpha_2\in \mathfrak g_2^*$, the $G_1$ action on $M_{\alpha_2}$ is multiplicity free;
        \item For every $\alpha_1\in \mathfrak g_1^*$, the $G_2$ action on $M_{\alpha_1}$ is multiplicity free.
    \end{enumerate}

    \begin{proof}
        By symmetry, it suffices to show that $(1)\Leftrightarrow (2)$. Let  
       $\alpha=\alpha_1+\alpha_2\in \mathfrak g^*$ and let $S$ be a piece of $M_{\alpha_2}$. For any $\alpha_1\in \mathfrak g_1^*$ we have that $S\cap \mu_{12}^{-1}(G_1(\alpha_1))/G_1$ is a stratified space with symplectic strata
        \begin{align} \label{pieces of doule reduction}
            \bigr(S\cap \mu_{12}^{-1}(G_1( \alpha_1))\bigr)/G_1=\bigcup_{i\in I} P_i.
        \end{align}
         Each of the $P_i$'s are pieces of $M_\alpha$.
         
         Let us now first assume that $M$ is a multiplicity free $G$-space. By Lemma \ref{lemma coisotrop iff red spac discrete}, we get that every $P_i$ consists of points. Again, using Lemma \ref{lemma coisotrop iff red spac discrete} we conclude that the $G_1$-action on every piece of $M_{\alpha_2}$ is multiplicity free. Now, if $f,g\in C^\infty(M_{\alpha_2})^{G_1}$, we have that the restriction of $\{f,g\}$ to any piece is identically zero. Thus, $\{f,g\}=0$ and $M_{\alpha_2}$ is a multiplicity free $G_1$-space.

         Conversely, assume that the $G_1$-action on $M_{\alpha_2}$ is multiplicity free for all $\alpha_2$. Write $\alpha=\alpha_1+ \alpha_2$ as above with $\alpha_1\in \mathfrak g^*_1$ arbitrary. By (\ref{pieces of doule reduction}), each piece of $M_\alpha$ is a symplectic stratum of the reduced space at $\alpha_1$ for the $G_1$-action on a piece $S$ of $M_{\alpha_2}$. By hypothesis, and using again Lemma \ref{lemma coisotrop iff red spac discrete}, we have that for each piece $S$ of $M_{\alpha_2}$, the reduced space $\bigr(S\cap \mu_{12}^{-1}(G_1(\alpha_1))\bigr)/G_1$ is discrete. Hence, all the pieces of $M_\alpha$ are discrete and $M$ is a multiplicity free $G$-space.
    \end{proof}
\end{thm}

\section{Hamiltonian actions on K\"ahler manifolds}\label{sect kahler}
Let $(M,\omega, J,\mathtt h)$ be a K\"ahler manifold and let $G$ be a compact Lie group acting by holomorphic isometries with $M/G$ connected. We assume that the $G$-action on $M$ extends to an holomorphic action of the complexification $G^\C$ on $M$ and the $G$-action on $M$ is Hamiltonian and so there exists a momentum map $\mu:M\lra \lieg^*$. We fix $\langle \cdot,\cdot \rangle$ an $\mathrm{Ad}(G)$-invariant scalar product on $\lieg$ and so we may think the momentum map as $\lieg$-valued map $\mu:M \lra \lieg$ by means of $\langle \cdot,\cdot \rangle$. The following result is well-known and it follows from the equvariance property of the momentum map.
\begin{lemma}\label{equi}
Let $x\in M$. Then $\mu(x)\in \lieg^{G_x}$. In particular, $[\mu(x),X]=0$ for any $X\in \lieg_x$.
\end{lemma}
\begin{defini}
A point $x\in M$ is called:
\begin{enumerate}
    \item \emph{stable} if $\lieg^\C_x=\{0\}$ and $G^{\C}(x) \cap \mu^{-1}(0) \neq \emptyset$;
    \item \emph{semistable} if $\overline{G^{\C}(x)} \cap \mu^{-1}(0) \neq \emptyset$;
    \item \emph{polystable} if $G^{\C}(x) \cap \mu^{-1}(0) \neq \emptyset$.
\end{enumerate}
\end{defini}
We denote by $M^{s}$, $M^{ss}$, $M^{ps}$ the set of the stable points, respectively semistable, polystbale points. These sets are obviously $G$-invariant. Therefore we talk about stable, semistable and polystable $G$-orbit. The results of Heinzner-Huckleberry-Loose \cite{h,hh,hhl,hl}, see also Sjamaar \cite{sja}, show that 
\begin{thm}\label{good-quotient}
The subsets $M^{s}$, $M^{ss}$ are open in $M$. The closure in $M^{ss}$ of a $G$-orbit contains a unique polystable orbit. Moreover there exists a good quotient $\pi:M^{ss} \lra Q$ with the properties:
\begin{enumerate}
    \item the induced map $\mu^{-1}(0)/G \lra Q$ is an homemorphism.
    \item two semistable $G$-orbits have the same image in $Q$ if and only if their closures in $M^{ss}$ are not disjoint, and this happens if and only if
the polystable orbits in their closures in $M^{ss}$ coincide.
\end{enumerate}
\end{thm} 
Therefore, the set of the semistable points admits a good quotient which can be identified as a topological space with the corresponding (possibly singular) symplectic quotient. 
There is no compactness or completeness
condition needed.

Let $H$ be a principal isotropy and let $F$ be a connected component of the core $\core$. Let $G_F$ denote the subgroup of $G$ which preserves $F$. Then $N(H)^o \subseteq G_F \subseteq N(H)$, $F$ is a K\"ahler manifold and the $G_F$-action on $F$ is by holomorphic isometries. Since $F$ is K\"ahler it follows that $F$ is $G_F^\C$-invariant. By Lemma \ref{equi}, keeping in mind $F\subseteq M^H$, we have $\mu(F)\subseteq \mathfrak{z}(H)$, where $\mathfrak{z} (H)$ denotes the Lie algebra of the centralizer of $H$. Since $Z(H) \subseteq N(H)$ we get the following result
\begin{prop}\label{mr}
The $G_F$-action on $F$ is Hamiltonian and a momentum map is given by the restriction of the momentum map $\mu$ to $F$.
\end{prop}
Summing up, the $G_F$-action on $F$ is Hamiltonian with a momentum map $\mu:F \lra \mathfrak{z}(H)$ and the $G_F$-action on $F$ extends to a holomorphic action of the complexification $G_F^\C$ on $F$. Altought $G_F$ is not connected, Theorem \ref{good-quotient} holds for the $G_F^\C$-action on $F$, see for instance \cite[Theorem 1.1]{PG}. Theorem \ref{reduction2} shows that inclusion
\[
\mu_{\vert_{F}}^{-1}(0)=\mu^{-1}(0) \cap F \hookrightarrow \mu^{-1}(0)
\]
induces a homeomorphism
\[
\mu^{-1}(0) \cap F / G_F \longrightarrow \mu^{-1}(0) /G.
\]
and so we get the following result.
\begin{thm}
The good quotient of the semistable points of $M$ is homemorphic to the set of semistable points of $F$.
\end{thm}
Let $x\in M$, let $\beta \in \mathfrak g$ and let $\mu^\beta:M \lra \R$ denote the contraction of the momentum along $\beta$, that is, $\mu^\beta (x):= \langle \mu(x),\beta \rangle$. Then $\mathrm{grad}\, \mu^\beta =J (\beta_M )$ and the function
\[
\lambda(x,\beta,t)=\langle \mu(\exp(t\textbf{i}\beta)x, \beta \rangle,
\]
increases, and so the limit
\[
\lambda(x,\beta)=\lim_{t\to+\infty} \lambda(x,\beta,t)
\]
exists.
\begin{defini}
The maximal weight of $x\in M$  in the direction of $\beta$ is the numerical
value
\[
\lambda(x,\beta)=\lim_{t\to+\infty} \lambda(x,\beta,t) \in \R \cup \{+\infty\}
\]    
\end{defini}
Let $x\in M$ and let $\beta \in \lieg$.  Let $c_\beta (t)=\exp (t \mathbf{i} \beta) x$. The energy of $c_\beta$ is given by 
\[
E(c_\beta)=\int_{0}^{+\infty} \mathtt h ( \mathrm{grad}\, \mu^\beta (c_\beta (t)), \mathrm{grad}\, \mu^\beta (c_\beta (t)) ) \mathrm{dt},
\]
where $\mathtt h$ denotes the K\"ahler metric.
The following definition was introduced by Teleman \cite{te}.
\begin{defini}
A $G$-action on $M$ is called \emph{energy complete} if for any $x\in M$ and for any $\beta\in \lieg$, if $E(c_\beta)<+\infty$ then $\lim_{t\mapsto +\infty} \exp(t\mathbf{i}\beta )x$ exists.    
\end{defini}
A Theorem of Teleman \cite[Theorem 3.3]{te}, see also \cite{bj,mu}, shows that
\begin{thm}\label{semistable}
If the $G$-action on $M$ is energy complete then $x\in M$ is semistable if and only if $\lambda(x,\beta)\geq 0$ for any $\beta \in \lieg$.
\end{thm}
Since $M^{ss}$ is $G^\C$-invariant it follows that $M^{ss}=G^\C (M^{ss} \cap F)$. 
\begin{thm}
If $M$ is energy complete then $M^{\mrm{ss}} \cap F$ is the set of semistable points of $F$ relative to the $G_F^\C$-action on $F$.
\end{thm}
\begin{proof}
Let $x\in F^{ss}$. Since the momentum map of the $G_F$-action on $F$ is the restriction of the momentum map of the $G$-action on $M$ it follows that $x\in M^{ss}$. 

Conversely, let $x\in M^{ss}\cap F$. By Theorem \ref{semistable} we have $\lambda(x,\beta)\geq 0$ for any $\beta \in \lieg$ and so for any $\beta \in \mathfrak n (\mathfrak h)$. Applying Theorem \ref{semistable}, we get $x\in F^{ss}$ concluding the proof.
\end{proof}
Let $\nu$ denote the norm square momentum map, that is, 
$\nu(p)=\frac{1}{2} \parallel \mu(p) \parallel^2$, where $\parallel \cdot \parallel$ is the norm relative to the $\mathrm{Ad}(G)$-invariant scalar product $\langle \cdot,\cdot \rangle$ on $\lieg$. $\nu$ restricted to $F$ is the norm square momentum map relative to the $G_F$-action on $F$ and will be denoted by $\nu_F$. Since $GF=M$ it follows that 
\[
\mathrm{inf}_{x\in M} \nu=\mathrm{inf}_{x\in F} \nu_F, \qquad 
\mathrm{sup}_{x\in M} \nu=\mathrm{sup}_{x\in F} \nu_F.
\]
By \cite[Corollary 6.12 page 178]{heinzner-schwarz-stoetzel} we get the following result
\begin{prop}
If $\nu$ has a local maximum at $x\in F$ then both $G(x)$ and $G_F (x)$ are complex.
\end{prop}
\begin{prop}
If $x\in F$ and $G(x)$ is complex then $G_F (x)$ is complex as well.   
\end{prop}
\begin{proof}
Since $G(x)=G^\C(x)$ and $G(x) \cap F=G_F (x)$, it follows that $G_F^\C (x)=G_F(x)$, concluding the proof.
\end{proof}
From now on, we assume that $M$ is compact.

Let $\mathcal C_M$ denote the set of critical points of $\mu$. By the $G$-invariance of the norm square gradient map the set of critical points $\mathcal C_M$ of $\nu$ is $G$-invariant and so $\mathcal C_M =G(\mathcal C_M \cap F)$. 
\begin{lemma}
$\mathcal C_M =G\mathcal C_F$, where $\mathcal C_F$ is the set of the critical points of $\nu_F$.
\end{lemma}
\begin{proof}
The gradient of the norm square gradient map at $x$ is given by  $\left(\mu(x)\right)_M (x)$. By Corollary \ref{mr} if $x\in F$ then $\left(\mu(x)\right)_M (x)\in T_p F$. This proves $\mathcal C_F=\mathcal C_M \cap F$ concluding the proof.
\end{proof}
Let $\mathcal B_M=\mu(\mathcal C_M)$, respectively $\mathcal B_F=\mu(\mathcal C_F)$ the set of the singular values of $\nu$, respectively the set of the singular values of $\nu_F$. By the above Lemma we get  $\mathcal B_M=G \mathcal B_F$. 

Let $\beta \in \mathcal B_M$ and Let $T_\beta :=\overline{\exp(\R \beta)}$. Since $N=\{x\in M:\, \beta_M (x)=0\}=M^{T_\beta}$ it follows that any connected components of $N$ is a totally geodesic submanifold of $M$. The function $\mu^\beta$ is constant on any connected component of $N$. 

Let $Z_\beta:=\{x\in M:\, \mu^\beta (x)=\parallel \beta \parallel^2\}$ and let $C_\beta :=G(Z_\beta \cap \mu^{-1}(\beta))$. Let $Y_\beta:=\{x\in M:\, \lim_{t\mapsto -\infty} \exp(t\textbf{i}) \beta)x\in Z_\beta\}$. The stratum $S_\beta$ associated to $\beta$ is a $G$-invariant neighborhood of $C_\beta$ in $GY_\beta$ \cite{kirwan}. A proof of the the following result is given in \cite{heinzner-schwarz-stoetzel,kirwan}.
\begin{prop}
Let $\beta_1,\beta_2 \in \mathcal B_M$. The following are equivalent:
 \begin{enumerate}
 \item $C_{\beta_1} = C_{\beta_2}$;
 \item $\beta_1 \in \mathrm{Ad}(G)(\beta_2)$;
 \item $S_{\beta_1} =S_{\beta_2}$.
 \end{enumerate}
\end{prop}
Hence, the stratum $S_\beta$ only depends on $\mathrm{Ad}(G)(\beta)$. Since $M$ is compact there are only finitely many strata, see for instance \cite{heinzner-schwarz-stoetzel,kirwan}.

By the above Proposition we may assume $\beta \in \mathcal B_F$. 

Let $x_0 \in M$ and let $x:\R \lra M$ be the negative gradient flow of $\nu$ through $x_0$. Since $\mathrm{grad}\, \nu (x)=\left(\mu(x)\right)_M (x)$ it follows that the complexified orbits are invariant under the gradient flow. 
Duistermaat, see Lerman \cite{lgf} and also \cite{ds}, noted that the Lojasiewicz gradient inequality holds for the norm square momentum map. In particular, the limit
\[
x_\infty :=lim_{t\to +\infty} x(t),
\]
exists and satisfies $\left(\mu(x_\infty)\right)_M (x_\infty)=0$. Moreover, the $G$-orbit through $x_\infty$ depends only of the $G^\C$-orbit through  $x_0$, see \cite[Theorem 6.4 page 43]{ds} and also \cite[Theorem 3.3. page 12]{bjc}.
\begin{prop}
Let $\beta \in \mathcal B_F$. Then $S_\beta \cap F$ is the stratum associated to $\beta$ relative to $\nu_F$.   
\end{prop}
\begin{proof}
Let $q\in S_\beta \cap F$. Since the momentum map of the $G_F$-action on $F$ is the restriction of the momentum map $\mu$ to $F$ it follows that negative gradient flow of $\nu$ coincides with the negative gradient flow of $\nu_F$. This implies  that $S_\beta \cap F$ is contained in the stratum associated to $\beta$  relative to $\nu_F$. 

Let $\tilde{S}_\beta$ be the stratum associated to $\beta$. By Corollary $7.6$ $p.24$ in \cite{heinzner-schwarz-stoetzel} we have
\[
\tilde{S}_\beta=\{p\in F:\, \beta \in \mu(\overline{G_F^\C \cdot p})\, \mathrm{and}\, \parallel \beta \parallel \leq \parallel \mu(gp)\parallel\, \forall g \in G_F^\C\}
\]
and so  $\tilde{S}_\beta \subseteq S_\beta \cap F$, concluding the proof.
\end{proof}
Summing up, we have proved the following result.
\begin{thm}
Let $\nu_F :F \lra \R$ be the norm square momentum map. Let $B_F=\mu(\mathcal C_F)$ be the critical value of $\nu_F$. Let $S_\beta^F$ denote the stratum of the Morse-like function $\nu_F$ associated to $\beta \in \mathcal B_F$. Then 
\[
M=\bigsqcup_{\beta \in B_F} G^\C S_\beta^F,
\]
is the smooth stratification of $M$ induced to $\nu$.
\end{thm}

\end{document}